\documentclass[preprint,10pt]{elsarticle}

\usepackage[utf8]{inputenc}
\usepackage{amsmath}
\usepackage{amssymb}
\usepackage{algorithmic}
\usepackage{graphicx}
\usepackage[ruled,vlined,linesnumbered]{algorithm2e}
\usepackage{appendix}
\usepackage{xcolor}

\newcounter{todocounter}

\SetAlFnt{\footnotesize}

\begin{document}

\begin{frontmatter}

\journal{}

\title{Approximating linear response of physical chaos}
\author{Adam A. \'Sliwiak\corref{cor}}
\ead{asliwiak@mit.edu}
\author{Qiqi Wang}
\ead{qiqi@mit.edu}
\address{Center for Computational Science and Engineering\\ Department of Aeronautics and Astronautics\\ Massachusetts Institute of Technology\\ 77 Massachusetts Avenue, Cambridge, MA, 02139, USA}
\cortext[cor]{Corresponding author.}

\begin{abstract}
    Parametric derivatives of statistics are highly desired quantities in prediction, design optimization and uncertainty quantification. In the presence of chaos, the rigorous computation of these quantities is certainly possible, but mathematically complicated and computationally expensive. Based on Ruelle's formalism, this paper shows that the sophisticated linear response algorithm can be dramatically simplified in higher-dimensional systems featuring a statistical homogeneity in the physical space. We argue that the contribution of the SRB (Sinai-Ruelle-Bowen) measure change, which is an integral part of the full linear response, can be completely neglected if the objective function is appropriately {\it aligned} with unstable manifolds. This abstract condition could potentially be satisfied by a vast family of real-world chaotic systems, regardless of the physical meaning and mathematical form of the objective function and perturbed parameter. We demonstrate several numerical examples that support these conclusions and that present the use and performance of a reduced linear response algorithm. In the numerical experiments, we consider physical models described by differential equations, including Lorenz 63, Lorenz 96, and Kuramoto-Sivashinsky.         
\end{abstract}
\begin{keyword}
Chaos, Ruelle's linear response theory, Sensitivity analysis, Space-split sensitivity (S3), SRB measure gradient
\end{keyword}
\end{frontmatter}

%Introduction----------------------------
\section{Introduction}\label{sec:intro}%To be filled out later
%Introduction to linear response
{\it Linear response theory} (LRT) \cite{book-lrt} provides an array of mathematical methods for analysis of system's reaction to small perturbations of imposed forces or control parameters. In particular, linear response of a dynamical system should be understood as the derivative of its output with respect to an input parameter. The name "linear response" is a direct consequence of the Taylor series expansion, which indicates that the system's reaction can be approximated by a linear function involving two terms: the unperturbed term and parametric derivative re-scaled by the imposed perturbation. Indeed, the use of Taylor series reveals one fundamental aspect of LRT. Namely, based only on information about the system in the unperturbed state, its response can be predicted for any small perturbation. Consequently, LRT is applicable to systems that vary differentiably with respect to its input. Efficient numerical algorithms for approximating linear response are fundamental in design optimization, uncertainty quantification, control engineering and inverse problems. These LRT-based computational tools are used in several fields of physics: electromagnetism \cite{kontani-electromagnetics}, plasma physics and fusion \cite{haskey-plasma}, statistical physics \cite{lucarini-statphysics}, turbulent flows \cite{larsson-turbulence}, climate dynamics \cite{ragone-climate}, and many more.  

%Linear response of chaos - Ruelle
In the presence of chaos, the classical formulation of LRT is modified. The reaction of a chaotic system is measured in terms of certain statistical quantities, e.g., long-time averages. Under the assumption of {\it ergodicity}, the statistics do not depend on initial conditions. Therefore, for a given chaotic model, the long-time statistics can be manipulated only by varying the input parameters. A prominent result in the field of LRT is the work of Ruelle \cite{ruelle-original,ruelle-corrections}, who rigorously derived a closed-form expression for the linear response of chaos. The major assumption of Ruelle's derivation is {\it uniform hyperbolicity}, which is a mathematical idealization of chaotic behavior. We postpone the description and explanation of this property for the following section of the paper. Solid numerical evidence found in literature clearly indicates that uniform hyperbolicity is a sufficient, bot not necessary, condition for the differentiability of statistics \cite{blonigan-phdthesis,chandramoorthy-phdthesis}. Indeed, these empirical results are consistent with {\it hyperbolic hypothesis} of Galavotti and Cohen \cite{galavotti-hypothesis}. That hypothesis presumes that several high-dimensional chaotic systems behave as though they were uniformly hyperbolic. It does not mean, however, that all properties of uniform hyperbolicity are satisfied by those systems, but several consequences following from this fundamental assumption could still be valid. This was clearly demonstrated in \cite{ni-jfm}, where the author argued that the long-time averages computed for a 3D turbulence model are smooth despite local non-hyperbolic behavior.          

%Different approaches of computing linear response - approximative 

While Ruelle's theory is regarded as one of the cornerstones in the field, its original expression for linear response is impractical due to the {\it butterfly effect}, i.e., exponential growth of tangent solutions in time. The ensemble method proposed in \cite{chandramoorthy-ensemble} circumvented this problem by computing ergodic averages along several truncated trajectories. Despite its simplicity, this method suffers from prohibitive computational costs induced by large variances of partial sensitivities. Shadowing methods \cite{wang-lssconvergence,ni-nilss} depart from the direct evaluation of Ruelle's expression by approximating the shadowing trajectories \cite{pilyugin-shadowing}, which lie in close proximity to the original path for a long period of time. Methods of this type have successfully been applied to high-dimensional fluid mechanics systems \cite{blonigan-phdthesis,ni-jfm}. However, a recent study \cite{chandramoorthy-shadowing} demonstrated that shadowing trajectories may be {\it nonphysical} and that their statistical behavior could be dramatically different than that of the reference trajectory. This unwanted behavior had also been observed in earlier studies, e.g. in \cite{blonigan-phdthesis}, which demonstrated large errors in shadowing-based approximations in spite of the apparently smooth behavior of the statistics. To the best of our knowledge, no rigorous studies that quantify or bound shadowing errors due to the problem of {\it nonphysicality} are available. An alternative way of computing linear response involves the Fluctuation-Dissipation Theorem (FDT) \cite{kubo-fdt}, which provides a time-convolution expression for the parametric derivative of statistics. FDT-based methods, such as the blended algorithm \cite{abramov-blended}, require some physics-informed assumptions to accurately reconstruct the linear response operator.

%Direct realizations of linear response - S3 and fast linear response
Recent algorithmic developments rely on the regularized variant of Ruelle's expression. Indeed, as originally proposed by Ruelle in \cite{ruelle-original}, one can apply integration by parts to the original formula in order to eliminate the product of Jacobians whose norm grows exponentially fast. However, since that formula involves Lebesgue integrals with respect to the Sinai-Ruelle-Bowen (SRB) measure \cite{young-srb} that is absolutely continuous only on unstable manifolds, an extra step is required before partial integration is applied. Namely, the input perturbation should be decomposed into two terms arranged in line with unstable and stable manifolds of the underlying dynamical system \cite{ruelle-original}. In case of flows (time-continuous systems), the center manifold should also be taken into account in the perturbation splitting \cite{ruelle-flows}. Based on this idea of regularization of Ruelle's closed-form expression, two conceptually similar methods for linear response emerged in the past two years. Those are the fast linear response algorithm \cite{ni-fast} and space-split sensitivity (S3) algorithm \cite{chandramoorthy-s3-new,sliwiak-s3}. The former additionally uses shadowing methods to approximate one of the terms resulting from the perturbation splitting. The S3 method, on the other hand, does not introduce any approximations except for the ergodic-averaging required for the evaluation of Lebesgue integrals inherited from the original formula. Indeed, the S3 method rigorously converges as a typical Monte Carlo procedure for any uniformly hyperbolic system \cite{chandramoorthy-s3-new}. Nevertheless, both methods can be summarized as follows. Split linear response into two terms (or three terms if considering a flow), such that one uses solutions of a regularized tangent equation (free of the butterfly effect), while the second term requires computing the divergence on unstable manifolds. The unstable divergence directly follows from the partial integration on the expansive tangent subspace. One of the by-products is the SRB density gradient representing the divergence of SRB measure. This quantity is obtained by differentiating the measure preservation law, which effectively requires solving a series of regularized second-order tangent equations \cite{sliwiak-densitygrad,ni-fast}. Differentiation of SRB measures, either explicit or implicit, is by far the most complicated and expensive part of both algorithms.

%What do we do in this paper? Highlight importance and novelty [IMPORTANT]
In this paper, we investigate if and under what circumstances the complex numerical procedures for linear response could be simplified. In particular, we attempt to answer the fundamental question about the significance of the SRB measure change. Rich numerical evidence found in the literate suggests that the computation of the SRB density gradient is not necessary to accurately approximate linear response in a number of popular physical systems. For example, the aforementioned shadowing methods, which in fact regularize the tangent equation and do not compute the curvature of unstable manifolds, have been proven successful in 3D turbulence models \cite{blonigan-phdthesis, ni-jfm}. Moreover, a recent theoretical study in \cite{ni-approx} concludes that if both the input perturbation and objective function follow the multivariate normal distribution, the effect of the measure change is expected to decay proportionally to $\sqrt{m/n}$, where $m$ is the number of positive Lyapunov exponents (LEs), while $n$ denotes the system's dimension. That work, however, does not provide any numerical examples. We show that the contribution of the unstable divergence could potentially be negligible if the objective function is specifically {\it aligned} with the unstable manifold. The meaning of {\it alignment} in this context is rigorously explained later in this work. Our numerical examples indicate that it is not uncommon that the SRB measure change is large and even has infinite variance, while its contribution to linear response might be negligible at the same time. This paradox may have huge implications for approximating sensitivities in large physical systems. The only obstacle is an additional requirement for the objective function, which typically has a concrete physical meaning. Our argument is based on the fact that a vast family of practicable systems are statistically homogeneous in physical space. They include popular models governing climate dynamics \cite{karimi-lorenz96}, turbulence \cite{larsson-turbulence}, population dynamics \cite{vaupel-population}, and several other phenomena. For such systems, we have freedom in representing any spatially-averaged objective function, which effectively increases the probability of its alignment with a tangent subspace.    

%Paper overview
The structure of this paper is the following. In Section \ref{sec:s3}, we thoroughly review the space-split sensitivity (S3) algorithm for linear response with an emphasis on potential difficulties. Subsequently, in Section \ref{sec:appproximating-unstable}, we explain the concept of alignment of the objective function and analyze its major implications in the context of the unstable contribution. A numerical experiment demonstrating a negligible effect of SRB measure change is presented. In Section \ref{sec:appproximating-hd}, we conjecture that the alignment constraint is not an obstacle for higher-dimensional systems with statistical homogeneity. Based on our analysis, we propose a reduced variant of the S3 method and apply it to approximate linear response of the Lorenz 96 and Kuramoto-Sivashinsky models. Section \ref{sec:conclusions} concludes this paper. \ref{sec:algorithm} and \ref{sec:derivatives} provide further technical details of S3: algorithm mechanics, implementation and cost analysis.

%Space-split sensitivity------------------
\section{Space-split sensitivity (S3) method for chaotic flows}\label{sec:s3}

The purpose of this section is twofold. First, we review the main results of the linear response theory, i.e., Ruelle's closed-form expression and its computable realization, known as the space-split sensitivity. Second, we present an extension of S3 to general hyperbolic flows and critically analyze its properties and major implications in the context of higher-dimensional systems. 

Throughout this paper, we consider a parameterized $n$-dimensional {\it ergodic} flow,
\begin{equation}
    \label{eqn:s3-flow}
    \frac{dx}{dt} = f(x;s),\;\;\;x(0)=x_0,
\end{equation}
with $m\geq 1$ positive Lyapunov exponents, where $s$ is a real-valued scalar parameter. The value of $m$ approximates the dimension of the unstable (expanding) subspace, while particular LE values indicate the rate of exponential expansion/contraction \cite{arnold-les}. Due to the assumed ergodicity, the statistical behavior of the system does not depend on the initial condition $x_0$.

For a given smooth objective function $J:M\to \mathbb{R}$, our ultimate goal is to approximate the parametric derivative of the long-time average of $J$, defined as
\begin{equation}
    \label{eqn:s3-sensitivity}
    \frac{d\langle J\rangle}{ds} := \frac{d}{ds}\lim_{T\to\infty}\frac{1}{T}\int_{0}^T J(x(t;s))\,dt, 
\end{equation}
where $M$ denotes the $n$-dimensional manifold defined by Eq. \ref{eqn:s3-flow}. We assume $J$ does not depend on $s$.

\subsection{Ruelle's formalism and S3}\label{sec:s3-derivation}
Under the assumption of {\it uniform hyperbolicity}, Ruelle derived a closed-form expression for linear response. Before we review the formula itself, we first focus on the assumption. A chaotic system is uniformly hyperbolic if its tangent space can be split into three {\it invariant} subspaces: unstable, stable and neutral. The first one and second one are spanned by expanding and contracting directions of the tangent space and they correspond to positive and negative LEs, respectively. These two subspaces respectively involve all tangent vectors that exponentially increase and decay in norm along a trajectory. In this paper, we focus on autonomous flows and thus the tangent space also involves a neutral subspace that is parallel to the flow vector $f$ and corresponds to the zero LE. In certain cases, a PDE-related dynamical system may involve more than one zero LE. For example, consider the Kuramoto-Sivashinsky equation with periodic boundary conditions. In this case, the neutral subspace is geometrically represented by a two-dimensional manifold (surface) that is tangent to $f$ and spatial derivative of the solution at every point on the attractor. The key aspect of hyperbolicity is that the three subspaces are clearly separated from each other, which means that the smallest angle between them is far from zero everywhere on the attractor. Hyperbolic systems are structurally stable and admit the SRB measure $\mu$ \cite{young-srb}, which contains the statistical description of the dynamics.

Assuming the system defined by Eq. \ref{eqn:s3-flow} is uniformly hyperbolic, Ruelle's linear response formula applies and can be expressed as follows \cite{ruelle-original,ruelle-corrections},
\begin{equation}
    \label{eqn:s3-ruelle}
    \frac{d\langle J\rangle}{ds} = \sum_{t=0}^{\infty}\int_M D(J\circ\varphi^t)\cdot\chi\,d\mu,
\end{equation}
where $g\circ h:=g(h)$, $\chi = \partial_s\varphi\circ\varphi^{-1}$, $\varphi^t = \varphi(\varphi^{t-1})$, $\varphi^0 (x) = x$, while $D$ denotes the gradient operator (first derivative) in phase space. The diffeomorphic map $\varphi:M\to M$ can be interpreted as a time integrator of Eq. \ref{eqn:s3-flow}. For example, using the second-order explicit Runge-Kutta method (midpoint rule) with step size $\Delta t$, $\varphi$ is related to $f$ through the following relation,
\begin{equation}
    \label{eqn:s3-time-int}
    x_{k+1} = \varphi(x_k) = x_k + \Delta t\,f(x_k + \frac{\Delta t}{2}\,f(x_k)).
\end{equation}
Since the system is assumed to be ergodic, the Lebesgue integral with respect to measure $\mu$ can be approximated as,
\begin{equation}
    \label{eqn:s3-ergodicity}
    \int_M h(x)\,d\mu = \lim_{T\to\infty}\frac{1}{T}\int_{0}^{T} h(x(t))\,dt \approx \frac{1}{N}\sum_{k=0}^{N-1} h(x_k)
\end{equation}
for any observable $h\in L^1(\mu)$ and a sufficiently large sample size $N$. Thus, the right-hand side (RHS) of Eq. \ref{eqn:s3-ruelle} could potentially be approximated by computing a sufficiently long trajectory, ergodic-averaging the integrand per Eq. \ref{eqn:s3-ergodicity}, and truncating the infinite series. However, note that
\begin{equation}
    \label{eqn:s3-ruelle-problem}
    D(J\circ\varphi^t)\cdot\chi = (DJ)_t\cdot (D\varphi)_{t-1}\,(D\varphi)_{t-2}...D\varphi\,\chi.
\end{equation}
$(DJ)_t$ denotes the phase-space gradient of $J$ evaluated $t$ time steps into the future. To facilitate the notation, we will drop the parentheses, i.e., $(DJ)_t:=DJ_t$. Therefore, unless $\chi$ is orthogonal to the unstable subspaces, the norm grows exponentially fast with $t$,
\begin{equation}
    \label{eqn:s3-ruelle-norm}
    \|D\varphi_{t-1}\,D\varphi_{t-2}...D\varphi\,\chi\|\sim \mathcal{O}(\exp(\lambda_1 t)), \;\;\;\lambda_1 > 0,
\end{equation}
which means the direct evaluation of the RHS of Eq. \ref{eqn:s3-ruelle} is computationally infeasible. The rate of exponential growth is determined by the leading LE denoted by $\lambda_1$. Indeed, due to the butterfly effect, the derivative of the composite function $J\circ\varphi^t$ is the most problematic aspect of Ruelle's original expression. Moreover, integration by parts is prohibited in this case, because one would also need to differentiate the SRB measure $\mu$ in the direction of $\chi$. In general, the measure is absolutely continuous only on the expanding subspace \cite{young-srb}. Therefore, integration by parts would be possible only if $\chi$ belongs to unstable manifolds everywhere in $M$, which is generally not the case.

Motivated by the work of Ruelle \cite{ruelle-original,ruelle-corrections}, the authors of \cite{chandramoorthy-s3-new,chandramoorthy-phdthesis} proposed a new method, called the space-split sensitivity (S3), which regularizes Ruelle's series for systems with one-dimensional unstable subspaces ($m=1$). Based on its extension to general hyperbolic maps in \cite{sliwiak-s3}, we derive and describe a space-split approach for chaotic flows with unstable manifolds of arbitrary dimension ($m\geq 1$). The main idea of S3, proposed in the aforementioned studies, is to decompose the perturbation vector $\chi$ into three terms,
\begin{equation}
    \label{eqn:s3-splitting}
    \chi = \chi^u + \chi^c + \chi^s = \left(\sum_{i=0}^m c^i\,q^i\right) + \left(c^0\,f\right) + \left(\chi - \sum_{i=0}^m c^i\,q^i - c^0\,f\right),
\end{equation}
such that $\chi^u$ and $\chi^c$ strictly belong to the unstable and neutral/center subspaces, respectively. In this splitting, $c^i$, $i=0,...,m$ are some scalars that are differentiable on the unstable subspace defined by a local orthonormal basis $q^i$, $i=1,...,m$. From now on, the superscript shall indicate the index of an array's component. This notation does not imply exponentiation, unless explicitly stated otherwise. There are two major benefits of the perturbation splitting defined by Eq. \ref{eqn:s3-splitting}:
\begin{itemize}
    \item the unstable part of the linear response (the one involving $\chi_u$) can now be integrated by parts, because it only involves directional derivatives along unstable subspaces,
    \item we can always find $c^i$, $i=0,...,m$ through orthogonal projection such that the stable part (the one involving $\chi_s$) of the linear response can be approximated by solving a regularized tangent equation that is bounded in norm.
\end{itemize}

We begin from exploring the second benefit of the splitting. Using the chain rule, one can rigorously show that the linear response defined by Ruelle's series equals the ergodic average of $DJ\cdot v$, where $v$ is a solution to the inhomogeneous tangent equation with $\chi$ as the source term. Thus, by replacing $\chi$ with $\chi^s$ in Eq. \ref{eqn:s3-ruelle}, we conclude that
\begin{equation}
    \label{eqn:s3-stable}
    \sum_{t=0}^{\infty}\int_M D(J\circ\varphi^t)\cdot\chi^s\,d\mu = \int_M DJ\cdot v\,d\mu,
\end{equation}
where
\begin{equation}
    \label{eqn:s3-inhomog}
    v_{k+1} = D\varphi_k\,v_k + \left(\chi_{k+1} - \sum_{i=0}^m c^i_{k+1}\,q^i_{k+1} - c^0_{k+1}\,f_{k+1}\right).
\end{equation}
The subscript notation indicates the time step, i.e., $f(x(k\Delta t)):=f_k$, assuming uniform time discretization. To solve Eq. \ref{eqn:s3-inhomog}, we need to project out the unstable component of $v$, otherwise its norm will grow exponentially in time at the rate proportional to the largest LE. Moreover, we should also project out the component tangent to the center manifold to eliminate the increase of sample variances, which we illustrate later in Section \ref{sec:s3-critical}. Therefore, we enforce $v$ to be orthogonal to the unstable-center subspace by imposing a set of $m+1$ constraints at every point on the manifold. Let $r_{k+1} = D\varphi_k\,v_k + \chi_{k+1}$ and, therefore,
\begin{equation}
    \label{eqn:s3-constraint1}
    (f_{k+1}\cdot f_{k+1})\,c_{k+1}^0 = f_{k+1}\cdot\left(r_{k+1} - \sum_{i=1}^m c_{k+1}^i\,q_{k+1}^i\right),    
\end{equation}
\begin{equation}
    \label{eqn:s3-constraint2}
    c_{k+1}^i = q_{k+1}^i\cdot\left(r_{k+1}-c_{k+1}^0\,f_{k+1}\right),\,\,\,i=1,...,m.    
\end{equation}
Eq. \ref{eqn:s3-constraint1}-\ref{eqn:s3-constraint2} define a linear system with $m+1$ equations and $m+1$ unknowns ($c^i$, $i=0,1,...,m$). The system's matrix involves the $m\times m$ identity block $I$, while its Shur complement can be expressed as follows,
\begin{equation}
    \label{eqn:s3-shur}
    S_{k+1} = I - \frac{Q_{k+1}^T f_{k+1}(Q_{k+1}^T f_{k+1})^T}{f_{k+1}\cdot f_{k+1}},
\end{equation}
where $Q$ is a an $n\times m$ matrix containing an orthonormal basis of the unstable manifold, $q^i$, $i=1,...,m$.
Thus, the coefficients $c^i$, $i=1,...,m$, stored in the array $c$ are obtained by solving the following reduced system,
\begin{equation}
    \label{eqn:s3-shur-system}
    S_{k+1}\,c_{k+1} = Q_{k+1}^T\,\left(r_{k+1} - \frac{f_{k+1}\cdot r_{k+1}}{f_{k+1}\cdot f_{k+1}}f_{k+1}\right),
\end{equation}
while $c^0$ is computed directly from Eq. \ref{eqn:s3-constraint1}. We conclude the stable part of the linear response can be evaluated through the ergodic average of $DJ\cdot v$ (see Eq. \ref{eqn:s3-ergodicity}), where $v$ satisfies Eq. \ref{eqn:s3-inhomog}--\ref{eqn:s3-constraint2}. 

The next step is the neutral contribution, which involves the perturbation component that is parallel to $f$. Analogously to Eq. \ref{eqn:s3-ruelle-problem}, we can expand,
\begin{equation}
    \label{eqn:s3-center-product}
    D(J\circ\varphi^t)\cdot\chi^c = D(J\circ\varphi^t)\cdot(c^0\,f) = c^0\,DJ_t\cdot \left(D\varphi_{t-1}...D\varphi\,f\right).
\end{equation}
Applying the Taylor series expansion, we note that
\begin{equation}
    \label{eqn:s3-center-indentities1}
    f(\varphi(x)) = f(x) + Df(x)\,(\varphi(x)-x) + \mathcal{O}((\varphi(x)-x)^2),
\end{equation}
and, analogously,
\begin{equation}
    \label{eqn:s3-center-indentities2}
    \varphi(x) = x + \Delta t\,Df(x) + \mathcal{O}(\Delta t^2).    
\end{equation}
By differentiating Eq. \ref{eqn:s3-center-indentities2} and plugging it to Eq. \ref{eqn:s3-center-indentities1}, we notice that in the limit $\Delta t\to 0$ we retrieve the covariance property, which reads
\begin{equation}
    \label{eqn:s3-center-covariance}
    f(\varphi(x)) = D\varphi(x)\,f(x).    
\end{equation}
This implies that the neutral part can be simplified to
\begin{equation}
    \label{eqn:s3-center-final}
   \sum_{t=0}^{\infty}\int_M D(J\circ\varphi^t)\cdot\chi^c\,d\mu = \sum_{t=0}^\infty \int_M c^0 DJ_t\cdot f_t\,d\mu = \sum_{t=0}^\infty \int_M c^0_{-t} DJ\cdot f\,d\mu.
\end{equation}
Eq. \ref{eqn:s3-center-final} means that the neutral part of the linear response equals the infinite series of $k$-time correlations between $c^0$, which is computed for the stable part, and $DJ\cdot f$. Under the assumption of uniform hyperbolicity, for any two H\"older-continuous observables $J$ and $h$, $k$-time correlations exponentially converge to the product of expected values as $t\to\infty$ \cite{chernov-correlations,young-srb}, i.e.,
\begin{equation}
    \label{eqn:s3-correlations}
    \left|\int_M (J\circ\varphi^t)\,h\,d\mu - \int_M J\,d\mu\,\int_M h\,d\mu\right|\leq C\delta^t
\end{equation}
for some $C>0$ and $\delta\in(0,1)$. In the context of the linear response theory, at least one of the observables has zero expectation with respect to $\mu$.
Using this property, we approximate the neutral part by truncating the infinite series and computing each Lebesgue integral through Eq. \ref{eqn:s3-ergodicity}. 

The final missing contribution of the total linear response is the unstable term. Indeed, this is the only term we can apply integration by parts to, which yields \cite{sliwiak-s3}
\begin{equation}
    \label{eqn:s3-unstable}
    \begin{split}
    &  \sum_{t=0}^{\infty}\int_M D(J\circ\varphi^t)\cdot\chi^u\,d\mu =  \sum_{t=0}^{\infty}\sum_{i=0}^m\int_M c^i \partial_{q^i}(J\circ\varphi^t)\,d\mu \\
    &
    =-\sum_{t=0}^{\infty}\sum_{i=1}^m\int_M (J\circ\varphi^t)\left(c^i\,g^i + b^{i,i}\right)\,d\mu,
    \end{split}
\end{equation}
where
\begin{equation}
    \label{eqn:s3-unstable-defs}
    b^{i,j}:=\partial_{q^j}c^i,\;\;\;g^i := \frac{\partial_{q^i}\rho}{\rho},
\end{equation}
the operator $\partial_{q^i}(\cdot):=D(\cdot)\cdot q^i$ denotes the directional derivative along $q^i$ in phase space, while $\rho$ denotes the density of the SRB measure $\mu$ conditioned on an unstable manifold. Several intermediate steps are required to derive the RHS of Eq. \ref{eqn:s3-unstable}. First, the SRB measure is disintegrated across parameterized unstable manifolds. Second, partial integration is applied within each parameterized subspace. The resulting boundary terms vanish as proven in \cite{ruelle-original}, which implies that in all integral transformations of this type, the boundary integrals can be neglected. The reader is also referred to \cite{sliwiak-srb} for a detailed description of every step of this process and relevant numerical examples. The major implication of Eq. \ref{eqn:s3-unstable} is that the composite function $J\circ\varphi^t$ is no longer differentiated, but there are two new quantities that must be computed instead. A rigorously convergent recursive algorithm for $b$ and $g$ has recently been proposed in \cite{sliwiak-s3}. That algorithm requires solving a collection of first- and second-order tangent equations, and was developed for discrete chaotic systems. In \ref{sec:algorithm}, we extend it to hyperbolic flows and analyze its cost. Notice that if $g$ and $b$ are available, then, analogously to the neutral part, the unstable term is expressed in terms of infinite series of $k$-time correlations. 

To summarize, the space-split method regularizes Ruelle's original expression by splitting it into three major parts: stable, neutral and unstable. Each of them can be approximated through ergodic-averaging of a single (in stable part) or many (in neutral and unstable parts) ingredients. Recent rigorous \cite{chandramoorthy-s3-new} and computational \cite{sliwiak-s3} studies have shown that the rate of convergence of all linear response parts is proportional to $1/\sqrt{N}$, where $N$ denotes the trajectory length. We highlight the fact that these studies were restricted to hyperbolic systems only. Thus, the S3 method is in fact a Monte Carlo procedure that relies on recursive formulas in the form of tangent equations that are executed to find $g$, $b$, $v$ and other necessary quantities.

\subsection{Numerical example: Lorenz 63}\label{sec:s3-example}
\begin{figure}
    \centering
    \includegraphics[width = 0.49\textwidth]{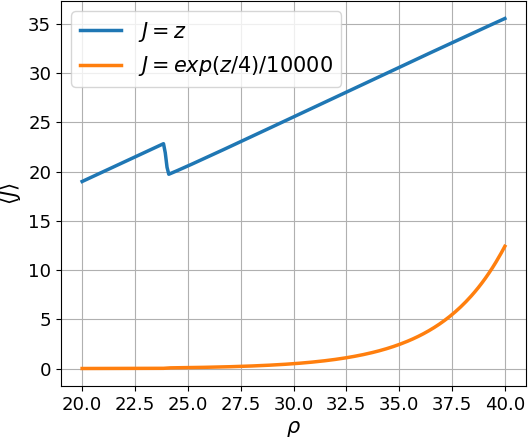}
    \includegraphics[width = 0.49\textwidth]{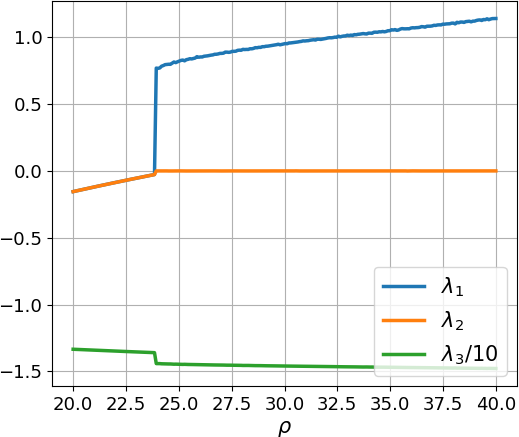}
    \caption{Long-time averages of two different objective functions (left) and Lyapunov exponents (right) versus the Rayleigh parameter $\rho$. Averages have been taken over $N\Delta t = 50,000,000$ and $N\Delta t = 5,000$ time units, respectively.}
    \label{fig:lorenz63_stats_les}
\end{figure}

%Lorenz 63 intro
To test the space-split algorithm (Algorithm \ref{alg:alg1}), we shall consider the three-dimensional Lorenz 63 system,
\begin{equation}
    \label{eqn:s3-lorenz63}
    \frac{dx}{dt} = \sigma(y-x),\;\;\;\frac{dy}{dt} = x(\rho - z) - y,\;\;\;\frac{dz}{dt} = xy - \beta z,
\end{equation}
which is one of the simplest chaotic flows. This ODE system models thermal convection of a fluid cell that is warmed from one side and cooled from the opposite side. The original study of this model \cite{lorenz-climate} demonstrated chaotic behavior at $\sigma=10$, $\beta = 8/3$, $\rho \gtrapprox 24$. For this choice of parameters, the strange attractor has a characteristic butterfly-shaped structure. The purpose of our experiment is to approximate the derivative of the long-time average of $J=J(z)$ with respect to the Rayleigh parameter $\rho$ using S3. In this section, $\rho$ should not be confused with the SRB measure density. Figure \ref{fig:lorenz63_stats_les} illustrates the behavior of the statistics of two different objective functions, as well as the three Lyapunov exponents for $\rho\in[20,40]$. We observe that $\lambda_1$ becomes positive for $\rho \gtrapprox 24$, which is consistent with the original study. The presence of a zero LE indicates there exists a tangent subspace that is parallel to the flow, which is typical for autonomous chaos. Note that, in the chaotic regime, both long-time averages seem to be differentiable in the considered parametric space.
%Sensitivity computation
To integrate Eq. \ref{eqn:s3-lorenz63} in time, we used the second-order explicit Runge-Kutta with step size $\Delta t = 0.005$. As described in \ref{sec:algorithm}, the space-split algorithm requires a few evaluations of first- and second-order differentiation operators of $\varphi$ every time step. For this particular time integrator, the computation of $D^2\varphi(\cdot,\cdot)$ involves three evaluations of the Hessian of $f$, per our derivations in \ref{sec:derivatives}. Fortunately, in case of the Lorenz 63 system, $D^2 f(\cdot,\cdot)$ is constant, which significantly reduces the cost. 

\begin{figure}
    \centering
    \includegraphics[width = 0.7\textwidth]{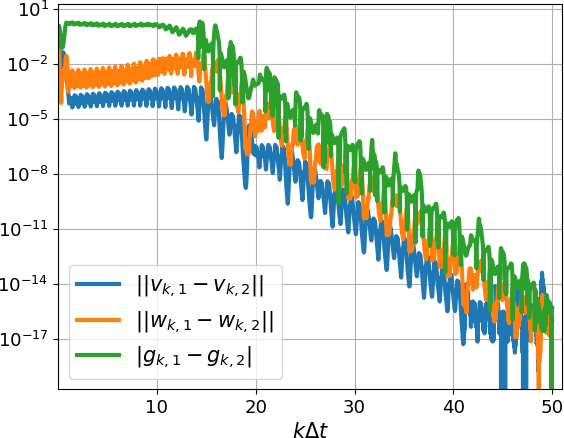}
    
    \vspace{1cm}\includegraphics[width = 0.7\textwidth]{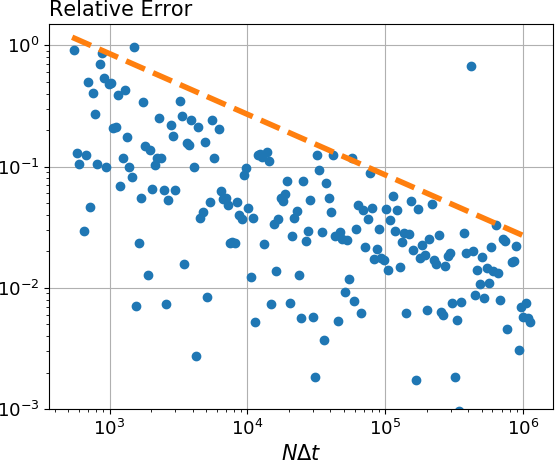}
    \caption{Upper: Relation of the norm/absolute value of quantities corresponding to two different initial conditions, labelled as 1 and 2, and time-averaging window $k\Delta t$. Lower: Relative error of the linear response approximation versus time-averaging window, computed for $J=z$ at $\rho = 28$. 200 independent simulations were run at a logarithmically uniform grid of $N\Delta t$. The dashed line represents a function $C/\sqrt{N\Delta t}$, $C>0$.}
    \label{fig:lorenz63_conv}
\end{figure}

The S3 algorithm relies on several recursive formulas in the form of recursive tangent equations. Earlier studies \cite{chandramoorthy-s3-new,sliwiak-s3} proved both analytically and numerically that these recursions converge exponentially fast in discrete hyperbolic systems. We numerically investigate if these results still apply the Lorenz 63 flow. The upper plot of Figure \ref{fig:lorenz63_conv} illustrates a convergence test for three different quantities: SRB density gradient $g$, tangent solution $v$ and its directional derivative (along $q$) $w$. These are three major ingredients that contribute to the total linear response. Along a single trajectory, we impose two different initial conditions for $v$, $w$ and $a$ (note $g=-q\cdot a$) and compute the norm/absolute value of the two solutions. The semi-logarithmic plot clearly indicates that all the norms decrease exponentially in time with a short transition at the beginning of simulation. To obtain a machine-precision approximation of these quantities, we need only 50 time units. Identical behavior have been observed for discrete systems \cite{sliwiak-s3}. We use this result to set the truncation parameter $T\Delta t = 100$ in our simulations to guarantee all ergodic-averaged quantities are very close to their true values. Another property of the S3 algorithm is the convergence rate of its final output, $\langle J\rangle/d\rho$, with respect the time-averaging window $N\Delta t$. Indeed, a truncation of the trajectory by choosing a finite $N$ is the only non-negligible source of error of the entire numerical procedure. The lower plot of Figure \ref{fig:lorenz63_conv} shows the decay of the relative error of the linear response approximation, which is computed with respect to the finite difference approximation of the slope of statistics generated in Figure \ref{fig:lorenz63_stats_les}. We observe that the error trend confirms theoretical predictions, which means that S3 behaves as a typical Monte Carlo simulation. 

\begin{figure}
    \centering
    \includegraphics[width = 0.7\textwidth]{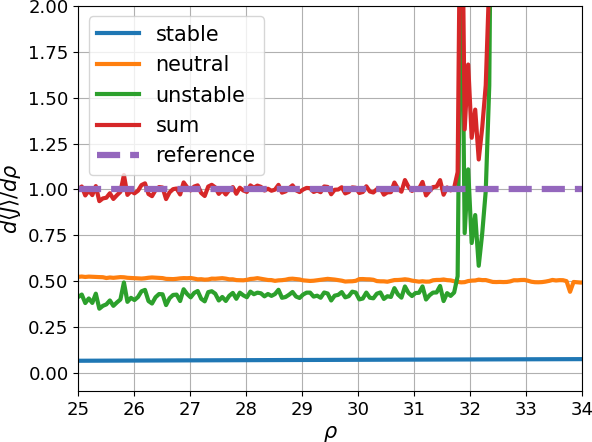}
    
    \vspace{1cm}\includegraphics[width = 0.7\textwidth]{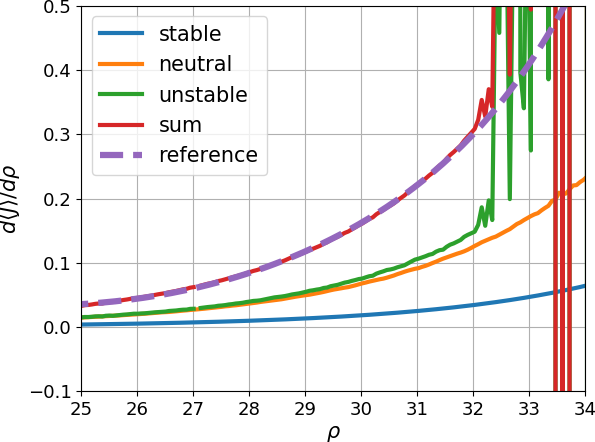}
    \caption{Output of Algorithm \ref{alg:alg1} generated for $J=z$ (upper) and $J=\exp(x/4)/10000$ (lower) at $144$ values of $\rho$ distributed uniformly. Each simulation was run for $N\Delta t = 1,000,000$ time units. The reference solution (dashed curve) was obtained using central finite differences and data shown in Figure \ref{fig:lorenz63_stats_les}. Before differentiation, we interpolated the data using first- and sixth-order polynomial fits, respectively.}
    \label{fig:lorenz63_sens}
\end{figure}

In our simulations, we truncate the infinite series by setting $K\Delta t = 50$, where $K$ represents the number of series terms contributing to the numerical approximation. The optimal value of $K\Delta t$ should be relatively small, given the exponential decay of correlations. In \cite{sliwiak-s3}, the reader will find a more detailed study about the impact of $K$ on the error. Based on the convergence study and our discussion above, we run Algorithm \ref{alg:alg1} for Lorenz 63 ($n=3$, $m=1$) to compute parametric derivatives of the long-time averages illustrated in Figure \ref{fig:lorenz63_stats_les} at $\rho\in[25,40]$. Figure \ref{fig:lorenz63_sens} shows the behavior of the obtained linear response approximations. For a wide range of Rayleigh constant values, S3 provides accurate estimations of the sensitivities. Indeed, for $\rho\in[25,32.3]$ we observe good agreement between the total sensitivity (denoted by ``sum") and corresponding reference values. At $\rho\approx 32.3$, the S3 approximation diverges due to the collapse of the unstable part. Note that, in both cases, the stable contribution is small compared to the two other terms. In the following section, we further explore the encountered problem and summarize critical aspects of the presented algorithm. 

\subsection{Critical view on S3}\label{sec:s3-critical}

In the context of approximating linear response of higher-dimensional chaos, we shall investigate potential problems of the S3 algorithm. In particular, we focus on dynamical properties of chaotic flows that might lead to numerical difficulties. Some algorithmic challenges, including the computational cost, are also discussed. 

\subsubsection{Special treatment of the neutral component}

\begin{figure}
    \centering
    \includegraphics[width = 0.6\textwidth]{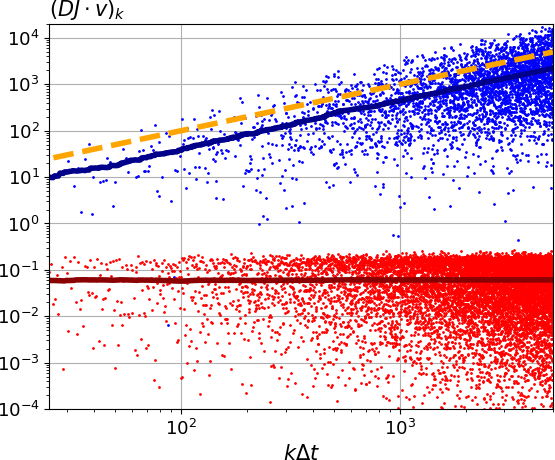}
    \caption{Discrete values of the stable integrand $DJ\cdot v$ computed using the S3 version described in Section \ref{sec:s3-derivation} (red) and its ``discrete" counterpart from \cite{sliwiak-s3} (blue). This simulation was performed for Lorenz 63 at $\rho = 28$. The solid lines represent the standard deviations of $DJ\cdot v$ collected from the beginning of the simulation until $k$th step. The dashed line represents a linear function.}
    \label{fig:disc-vs-cont}
\end{figure}

In Section \ref{sec:s3-derivation}, we derived a numerical scheme based on the three-term linear splitting in Eq. \ref{eqn:s3-splitting}. Indeed, there is a subtle difference between this splitting and the one proposed for discrete systems. In the former, the neutral term is treated separately thanks to which the stable term includes only tangent solutions that are parallel to the unstable-center subspace. In Figure \ref{fig:disc-vs-cont}, we plot discrete values of the stable integrand $DJ\cdot v$ obtained for Lorenz 63 at $\rho = 28$ using both versions of S3. We notice that if the neutral direction is not projected out from the tangent solution, then the standard deviation of $DJ\cdot v$ grows linearly with time. The extra projection against $f$ guarantees the standard deviation is approximately constant.

While the convergence of the Monte Carlo procedure is now guaranteed, the extra projection requires assembling, inverting, and differentiating the Schur complement. As described in \ref{sec:algorithm}, that minor conceptual adjustment requires major modifications of the ``discrete" version of S3.

\subsubsection{Problem with hyperbolicity and SRB measure gradient}

Recall that the major assumption of Ruelle's formalism is hyperbolicity. Any form of linearly separated perturbation splitting that enables partial integration and that guarantees boundedness of the stable part, e.g., the one presented in this paper or the shadowing-based variant proposed in \cite{ni-fast}, would be sufficient to construct stable numerical schemes. However, the dynamic structure of many chaotic flows, including the simple Lorenz 63 system, does not satisfy all basic properties of hyperbolicity. 

\begin{figure}
    \centering
    \includegraphics[width = 0.47\textwidth]{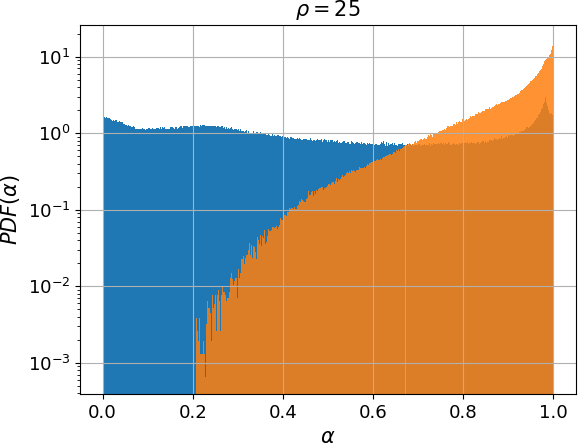}\hspace{5mm}
    \includegraphics[width = 0.47\textwidth]{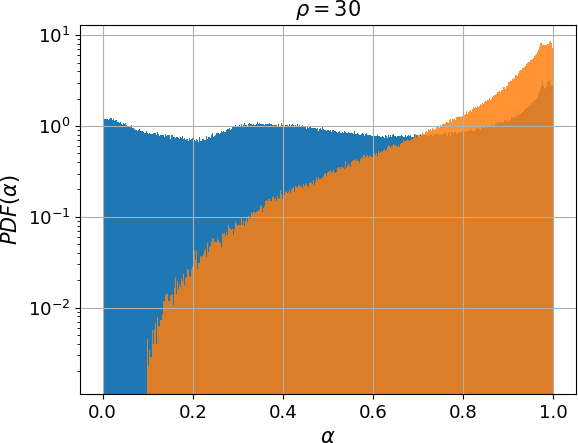}
    
    \vspace{1cm}
    \includegraphics[width = 0.47\textwidth]{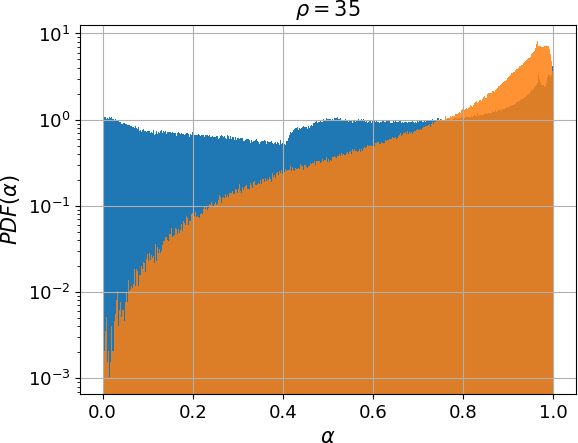}\hspace{5mm}
    \includegraphics[width = 0.47\textwidth]{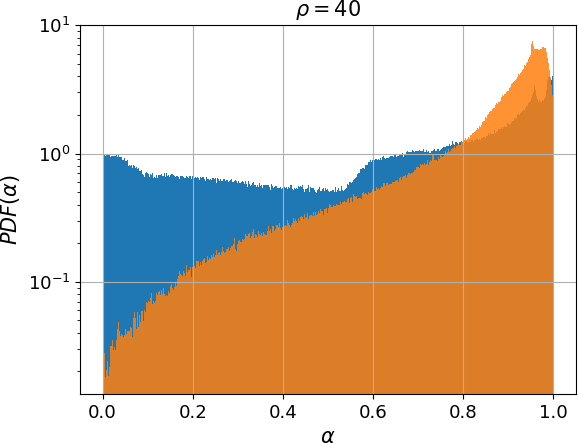}

    \caption{Distribution of the normalized measure $\alpha$ between unstable/center subspaces (blue PDF) and unstable-center/stable subspaces (orange PDF). They have been computed along a random trajectory of the Lorenz 63 system for 5000 time units. To increase the accuracy of PDFs, we used the fourth-order Runge-Kutta time integrator with $\Delta t = 0.005$.}
    \label{fig:hyperbolicity}
\end{figure}

In Figure \ref{fig:hyperbolicity}, we illustrate the distribution of tangency measures $0\leq\alpha\leq 1$ between two pairs of subspaces: 1) unstable and center, 2) unstable-center and stable, along a random trajectory of Lorenz 63 at different values of the Rayleigh parameter. To generate these plots, we used the fast algorithm for hyperbolicty verification proposed by Kuptsov in \cite{kuptsov-hyperbolicity}. The two measures we compute respectively represent $d_1$, and $2\,d_2$, which are rigorously defined by Eq. 7 in that work. The parameter $\alpha$ is closely related to the minimum angle between two subspaces normalized by $\pi/2$ as pointed out and tested in \cite{takeuchi-hyperbolicity}. If the statistical distribution of $\alpha$ is not strictly separated from the origin, i.e., it is very close to $\alpha = 0$, then several tangencies of a given subspace pair are highly likely to occur and vice versa. We observe that, regardless of the choice of $\rho$, there exist tangencies between the unstable and center subspaces. Several numerical examples presented in \cite{kuptsov-hyperbolicity} imply that the absence of unstable-center separation is a common property of several physical systems. However, for some $\rho$, the Lorenz 63 system admits splitting of the tangent space into unstable-center and stable subspaces. This behavior has been known in literature \cite{morales-singular} under the name of {\it singular hyperbolicity}. Note the Lorenz 63 loses this property at $\rho$ between $30$ and $35$. This coincides with the collapse of the S3 algorithm. In particular, the unstable term blows-up within this parameter regime, which suggests $\mu$ becomes rough along expansive directions within this interval. From the study on differentiability of statistics of the Lorenz system \cite{sliwiak-differentiability}, we learn that the SRB density gradient $g$ is Lebesgue integrable, i.e., $g\in L^1(\mu)$, if $\rho < 32$. If $\rho$ is close to the value of 28, then $g$ is even square-integrable. The authors of the same paper argue that, based on the earlier version of the S3 formulation and several numerical experiments, the integrability of $g$ implies differentiability of statistics and vice versa. We conclude that even if Eq. \ref{eqn:s3-ruelle} holds, one still needs to handle the by-products of partial integration, which poses a serious challenge for Monte Carlo algorithms such as S3.

The smoothness of the SRB measure is not guaranteed in non-hyperbolic systems, which means that some components of $g$ might not exist at all at some points on the attractor. Indeed, numerical experiments presented in \cite{kuptsov-hyperbolicity, takeuchi-hyperbolicity} indicate that some higher-dimensional physical systems, e.g., the Ginzburg-Landau equation, are clearly non-hyperbolic. Similar numerical results were provided for a 3D turbulent flow in \cite{ni-jfm}. Since $g$ is an integral part of the S3 procedure and its value is computed everywhere along a random trajectory, we expect that the unstable contribution might blow-up in case of such systems.

\subsubsection{Implementation and cost}
We shall now comment on practical aspects of the full linear response algorithm, which is described in \ref{sec:algorithm}. In terms of the implementation, both the stable and neutral parts do not require significant changes of the existing tangent/adjoint solvers. The former is obtained by solving a collection of first-order tangent equations. They are stabilized by step-by-step elimination of unstable-center tangent components through QR factorization that is needed to find a new basis of the subspace (matrix $Q$) and the Jacobian of coordinate transformation (matrix $R$). The $R$ factor can also be used to approximate $m$ largest LEs, which is indeed a very useful by-product of the algorithm. The unstable contribution requires the implementation of the second-order derivative operator, which is necessary for $g$ and $b$. While this is generally not a problem for simple systems, the need for a second-order tangent solver might require extra tools, such as automatic differentiation packages, for complicated higher-dimensional models. 

It turns out that the presence of the Hessian is not the major burden of the full S3 algorithm. The typical structure of large physical systems is sparse due to the localized stencils of the most popular spatial discretization schemes. Therefore, the computational cost of matrix-vector or tensor-vector products is typically linear in $n$. Two other factors that determine the total cost is the trajectory length $N$ and the number of positive LEs $m$. The former defines the accuracy of ergodic-averaging and indicates the number of primal/tangent solution updates, and thus contributes linearly to the total cost. Based on our estimate in \ref{sec:algorithm}, the final cost is proportional to the third power of $m$. The most expensive chunk of the algorithm is associated with the SRB density gradient $g$, which requires solving $\mathcal{O}(m^2)$ second-order tangent equations that is followed by a stabilizing normalization procedure consuming extra $\mathcal{O}(n\,m^3)$ flops. This might pose a serious challenge for systems with hundreds of unstable modes, such as 3D turbulence models.   

\subsubsection{Future prospects}

The non-approximative methods for computing linear response of chaotic systems, such as the S3 algorithm, provide a rich collection of numerical tools for analysis of the underlying dynamics. Its major drawback is that the derivation of its components relies on the assumption of hyperbolicity and smooth SRB measure. These properties might be violated leading to the collapse of some parts of the full S3 algorithm. Nevertheless, we acknowledge the growing popularity and interest in hyperbolic systems among physicists and engineers. In a comprehensive review book of Kuznetsov \cite{kuznetsov-hyperbolicity}, the author justifies this trend and provides several examples of hyperbolic attractors describing physical phenomena.

Despite the problems with hyperbolicty and large costs, can we still use some parts of the S3 algorithm to find accurate estimates of linear response for higher-dimensional systems? Figure \ref{fig:lorenz63_sens} indicates that both the neutral and stable contributions of Lorenz 63 remain ``stable" over the entire parametric regime. Moreover, as shown in \cite{sliwiak-differentiability}, the collapse of the recursion for $g$ does not mean Ruelle's linear response expression does no longer apply. Removal of the unstable contribution would dramatically reduce the cost of S3, as the expensive and potentially incomputable $g$ would no longer be needed. In case of Lorenz 63, however, the unstable contribution accounts for approximately $50\%$ of total sensitivity. Therefore, omission of the unstable contribution of this system would give rise to significant errors. This observation leads to a fundamental question. Are there systems whose unstable contribution is small and can be neglected? If so, are they relevant for practitioners? We try to answer those questions in the remainder of this paper.

%-------------Unstable contribution----------------------
\section{Unstable contribution: can we neglect it?}\label{sec:appproximating-unstable}

As we pointed out in Section \ref{sec:s3-critical}, the computation of the unstable part of linear response might be cumbersome due to several reasons. The purpose of this section is to provoke a discussion about the significance of that term. In particular, we shall present some evidence indicating that the unstable term could be negligible and thus completely neglected if certain conditions are met. 

Let us consider the leading term of Eq. \ref{eqn:s3-unstable}, i.e., the one corresponding to $t=0$,
\begin{equation}
    \label{eqn:approx-unstable-leading}
    U := \int_M J \, d\,d\mu,\;\;\;d:= d_{cg} + d_b :=\sum_{i=1}^m c^i\,g^i + b^{i,i}.
\end{equation}
Assuming the exponential decay of correlations holds, it is clear that the whole infinite series is small if $U$ is small. Applying the Cauchy-Schwarz and triangle inequalities, we upperbound the magnitude of $U$,  
\begin{equation}
    \label{eqn:approx-unstable-leading-estimate}
    |U| \leq  \|J\|_2\,\|d\|_2 \leq \|J\|_2\,\left(\|d_{cg}\|_2 + \|d_b\|_2\right),
\end{equation}
where $\|\cdot\|_2$ denotes the $L^2$ norm with respect to $\mu$ defined as
\begin{equation}
    \label{eqn:approx-l2norm}
    \|h\|_2 := \sqrt{\int_M h^2\,d\mu}
\end{equation}
for any scalar function $h\in L^2 (\mu)$. According to Inequality \ref{eqn:approx-unstable-leading-estimate}, we see that a small $L^2$ norm of the unstable divergence $d$ implies that the entire unstable contribution is negligible as well. Recall that the vector $c$ represents projections of the tangent solution $v$ onto the unstable subspace, which depends on both $\chi$, i.e., the parametric perturbation of the system, and geometry of the unstable manifold. The final term contributing to $\|d\|_2$ is the SRB density gradient, which represents measure change in $m$ orthogonal directions of the unstable subspace. These directions, stored in the $Q$ matrix, indicate how the unperturbed trajectory deforms in time. The rate of geometric expansion in the $i$-th direction is reflected by the $i$-th Lyapunov exponent $\lambda_i$, whose value can expressed in terms of the following ergodic average \cite{ershov-lyapunov},
\begin{equation}
    \label{eqn:approx-le-def}
    \lambda_{i} = \int_M \log\left|q^i(\varphi(x))\cdot D\varphi(x)\,q^i(x)\right|\,d\mu.
\end{equation}
We also acknowledge that the computation of $Q$ is an integral part of the S3 procedure (see \ref{sec:algorithm}). In that algorithm, the columns of $Q$ are sorted from the most expansive ($i=1$) to the least expansive ($i=m$) direction. Eq. \ref{eqn:approx-le-def} suggests that a bunch of infinitesimally close points will scatter very fast along the $q^i$ direction if $\lambda_i$ is large resulting in a small local measure change. In other words, larger expansion rates lead to the dilution of measure, which consequently decreases measure gradient. 
Therefore, assuming the positive LEs are separated from each other, we conjecture that the measure change along $q^1$ and $q^m$ are expected to be the smallest and largest, respectively. In particular, if $\lambda_1 > \lambda_2  ... >\lambda_m$, then
\begin{equation*}
    %\label{eqn:approx-unstable-gs}
    \|g^1\|_2 < \|g^2\|_2 < ... < \|g^m\|_2.
\end{equation*}
We verify this presumption later in a numerical experiment. Its major consequence is that we can potentially find two different directions on unstable manifolds along which the rates of change of $\mu$ are significantly different.      

As a side note, we bring up the fact that the two unstable contributions, associated with $d_{cg}$ and $d_{b}$, are the same in magnitude if $J\equiv1$. Indeed, using the definition of $b$, we observe that $\sum_{i=1}^m b^{i,i}:=\nabla_{\xi}\cdot c$, where $\nabla_{\xi}$ denotes the nabla operator (gradient) on unstable subspace. Thus, we can use Green's first identity to rewrite the latter term to
\begin{equation}
    \label{eqn:approx-unstable-b}
    \int_M J\,d_{b}\,d\mu = \int_M J\,\nabla_{\xi}\cdot c\,d\mu = -\int_M c\cdot \frac{\nabla_{\xi}(\rho J)}{\rho} d\mu \overset{J\equiv 1}{=} -\int c\cdot g\,d\mu,
\end{equation}
where $\rho$ denotes the measure density conditioned on a local unstable manifold. It is now evident that the two ingredients of $U$, $d_{cg}$ and $d_b$, involve both the array of $v$--$Q$ projections and a vector representing local relative measure change. The only difference between them is that, in the latter term, measure change is weighted by the value of $J$. If $J$ is not strongly-oscillatory nor has it large gradients in phase space, then $\nabla_\xi (\rho J)/\rho$ has a behavior similar to its non-weighted counterpart $g$.   

This analysis shows that there are two possible ways of reducing the norm of $U$, through $c$ and/or $g$. According to the definition of $c$, reducing its norm would restrict our analysis only to a certain parameter. Note that $c$ directly depends on $\chi$, which represents the parametric perturbation of the trajectory. On the other hand, $g$ contains information on statistics of the system in the unperturbed state. Therefore, neutralization of the effect of $g$ might allow us to dramatically decrease $|U|$, regardless of the parameter with respect to which the linear response is computed. In the remainder of this section, the concept of ``neutralization" will be explained in more detail. 

%Concept of alignment
Let us now consider a well-behaved objective function $J:M\to\mathbb{R}$, where $M$ is an orientable compact manifold. Let the tangent bundle of $M$ be expansive in all possible directions, which implies that all LEs are positive. Without loss of generality, we assume the volume integral of $J$ over $M$ is zero. Notice we can always add a constant number to $J$ to ensure the zero mean condition, as the constant shift does not affect linear response. Thus, $J$ can be expressed in terms of the divergence of a vector field $Z$, i.e.,
\begin{equation}
    \label{eqn:approx-divergence}
    J = \nabla_{\xi}\cdot Z.
\end{equation}
After plugging Eq. \ref{eqn:approx-divergence} to the expression for $U$, we can apply Green's first identity analogously to Eq. \ref{eqn:approx-unstable-b}, which yields
\begin{equation}
    \label{eqn:approx-divergence-green}
    U = - \int_M Z\cdot \frac{\nabla_\xi(\rho\,d)}{\rho}\,d\mu.
\end{equation}
Note that Eq. \ref{eqn:approx-divergence-green} contains all combinations of mixed second derivatives of the SRB measure. To minimize the effect of the measure change, we want to eliminate possibly as many components of $g$ as possible, especially those corresponding to the least expansive directions (highest indices). In an ideal scenario, we also want to neutralize the effect of those components of $g$ that remain. This could be achieved by choosing a $J$ that is {\it aligned} with $q^1$, which means that the statistics of $\nabla_\xi\cdot Z = \sum_{i=1}^m\,\partial_{q^i}\,Z^i$ is dominated by its first term ($i=1$), i.e., 
\begin{equation*}
    \|\partial_{q^1}Z^1\|_2 \gg \|\partial_{q^i}Z^i\|_2,\;i=2,...,m.
\end{equation*}
In this special case, we could approximate $U$ by keeping only the first term of $\nabla\cdot Z$. For the truncated expression, we apply integration by parts, which yields
\begin{equation}
    \label{eqn:approx-unstable-parts}
    U \approx U^1 := \int_M \partial_{q^1}Z^1\,d\,d\mu = -\int_M Z^1\,\left(d\,g^1 + \partial_{q^1} d\right)\,d\mu.
\end{equation}
The first benefit of the alignment is that we automatically eliminate second differentiation with respect to directions indicated by $q^2,...,q^m$ that correspond to the largest slopes of $\mu$. Therefore, the leading term of the unstable contribution is upperbounded as follows,
\begin{equation}
    \label{eqn:approx-unstable-leading-estimate2}
    |U^1| \leq  \|Z^1\|_2\,\left(\|d\,g^1\|_2\ + \|\partial_{q^1} d\|_2\right).
\end{equation}
The first term of the new inequality is proportional to $\|d\,g^1\|_2$. If $\|g^1\|_{\infty} \ll 1$, which is true if the measure is almost constant along $q^1$, then $\|d\,g^1\|_2 \ll \|d\|_2$. This scenario is very likely in systems with a broad Lyapunov spectrum. In the second term of Ineq. \ref{eqn:approx-unstable-leading-estimate2}, $d$ is differentiated in the most expansive direction $q^1$. It means that all components of the SRB density gradient weighted by $c$, are differentiated once more. This time, however, we differentiate in the direction of the mildest descent/ascent of $\mu$. One could visualize this process by considering the lateral boundary of a cylindrical solid. In this case, the tangent line computed along the solid's height is always parallel to the solid (zero slope). Any other slope is larger than zero. Differentiation of the non-zero slopes along the solid's height effectively kills them all. We can apply this analogy to our case, in which we differentiate one more in the direction of the smallest slope.  Therefore, the effect of the largest components of $g$ corresponding to the least expansive directions could be neutralized, in which case $\|\partial_{q^1} d\|_2$ is expected to be negligible. 

Through the above analysis, we conjecture that if $J$ is aligned with the most expansive direction of the unstable manifold, as defined above, and the positive part of the Lyapunov spectrum is not clustered around a certain value, it is possible to significantly reduce the magnitude of the unstable contribution. While the second condition is satisfied by many physical systems, the specific requirement for the objective function might be very restrictive. We now present a numerical example illustrating our argument.

%Numerical example

In our investigation, we will focus on the following $n$-dimensional chaotic map $\varphi:[0,2\pi]^n\to[0,2\pi]^n$ defined as
\begin{equation}
    \label{eqn:approx-sawtooth}
    x^i_{k+1} = 2\,x^i_k + s\,\sin(x^{i+1}_k - x^{i}_k) + t\,\sin(x^i_k)\;\mathrm{mod}\;2\pi,\,\,\,i=1,..,n,
\end{equation}
where $n\in\mathbb{Z}^+$, $s\in \mathbb{R}$, $t\in \mathbb{R}$ and $x^{n+1} = x^{1}$. This is an extension of the one-dimensional {\it sawtooth map} \cite{sliwiak-1d}, and therefore we shall refer to $\varphi$ defined by Eq. \ref{eqn:approx-sawtooth} as the {\it coupled sawtooth map}. The first term on the RHS introduces constant expansion that does not involve any parameters. Thus, if we set the coupling parameter to zero ($s=0$), we obtain $n$ independent maps with the same statistical behavior. If both the coupling and distorting terms are small, i.e., respectively $s$ and $t$ are small, then all Lyapunov exponents are clustered around the value of $\log 2$, which means that the attractor is expansive in all directions ($m=n$). By increasing $|s|$, we strengthen the coupling between the neighboring degrees of freedom. For $n=2$, the phase space gradient of the coupling term is parallel to the diagonal of the square manifold $[0,2\pi]^2$. Thus, the larger $|s|$, the stronger variations of the measure are expected along $[1,-1]^T$. In case of a weak distortion, i.e., $t\approx 0$, the SRB measure is expected to be approximately constant in the direction parallel to $[1,1]^T$. 

\begin{figure}
    \centering
    \begin{minipage}{1.0\textwidth}
        Case 1: $[s,t]=[0.05,0]$\\
        \includegraphics[width = 0.41\textwidth]{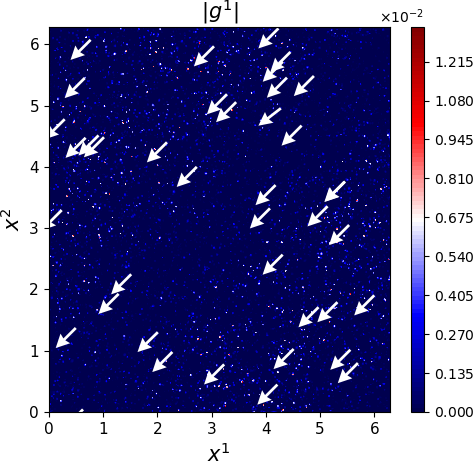}
        \hspace{4mm}
        \includegraphics[width = 0.4\textwidth]{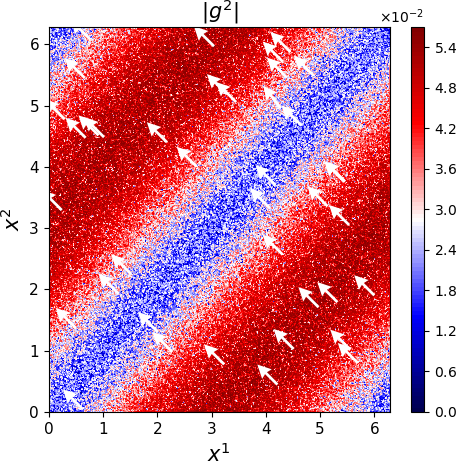}
    \end{minipage}%
    \vspace{2mm}
    \begin{minipage}{1.0\textwidth}
        Case 2: $[s,t]=[-0.75,0]$\\
        \includegraphics[width = 0.4\textwidth]{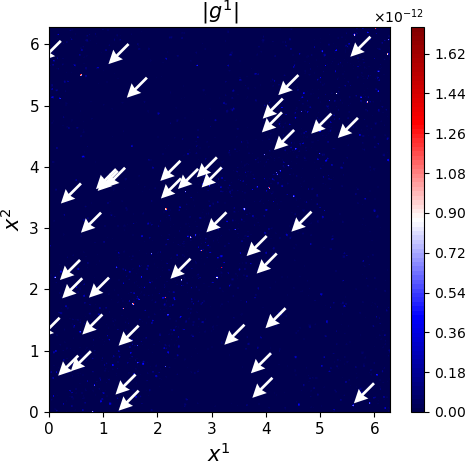}
        \hspace{4mm}
        \includegraphics[width = 0.4\textwidth]{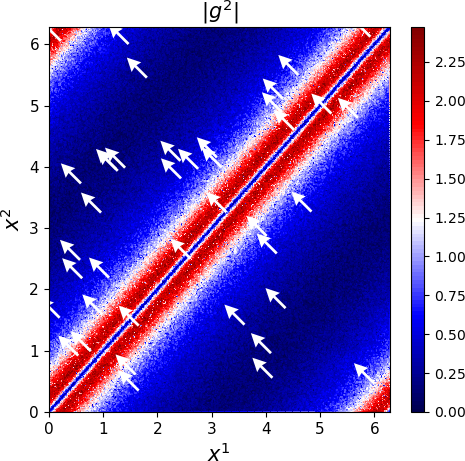}
    \end{minipage}%
    \vspace{2mm}
    \begin{minipage}{1.0\textwidth}
        Case 3: $[s,t]=[-0.75,0.5]$\\
        \includegraphics[width = 0.4\textwidth]{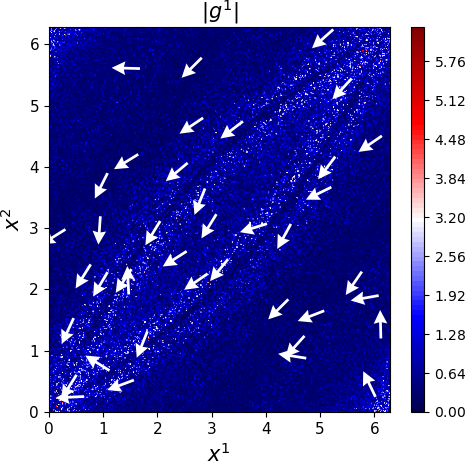}
        \hspace{4mm}
        \includegraphics[width = 0.4\textwidth]{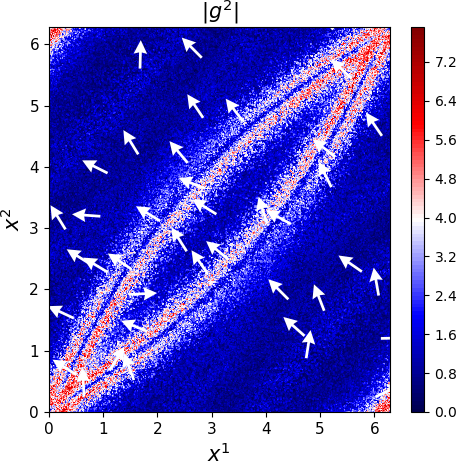}
    \end{minipage}
    
    \caption{Magnitude of two components of the SRB density gradient $g$ of the two-dimensional coupled sawtooth map with two positive LEs. White arrows respectively represent $q^1$ and $q^2$, which indicate local directions of differentiation. They are plotted every 5000 time steps. For each case, a trajectory of length $N = 3\cdot 10^5$ was generated.}
    \label{fig:cs-g}
\end{figure}

To verify these suppositions, we directly compute $g$ for $n=2$ at three different parameter sets: 1) $[s,t]=[0.05,0]$ (weak coupling, no distortion), 2) $[s,t]=[-0.75,0]$ (strong coupling, no distortion), 3) $[s,t]=[-0.75,0.5]$ (strong coupling combined with distortion). For this purpose, we use a part of the full S3 algorithm to compute $g$ along a trajectory (Lines 12-19 of Algorithm \ref{alg:alg1} in \ref{sec:algorithm}) and plot both $|g^1|$ and $|g^2|$ on $[0,2\pi]^2$. These results are illustrated in Figure \ref{fig:cs-g}. In all three cases, the first component of $g$ is statistically smaller in magnitude and features milder variations compared to the second one. They also confirm that the larger component of the relative measure change is approximately parallel to $[1,-1]^T$. Even in the presence of the distortion term (Case 3), the majority of white arrows, which indicate local directions $q^1$ and $q^2$, tend to be oriented diagonalwise. Notice that the larger coupling $|s|$, the larger rate of measure change in the least expansive direction represented by $q^2$. If there is no distortion and coupling is significant (Case 2), then the first component of $g$ is approximately zero everywhere in phase space. The largest measure gradients appear to be located around the $[1,1]^T$ diagonal. Furthermore, if the coupling weakens, then the rates of expansion along $q^1$ and $q^2$ become similar. In Case 1, the distribution of $g^1$ has geometric features similar to its counterpart. This is consistent with our analysis suggesting that both distributions are expected to have the same limits as $|s|\to 0$.   

\begin{figure}
    \centering
    \begin{minipage}{1.0\textwidth}
        \includegraphics[width = 0.45\textwidth]{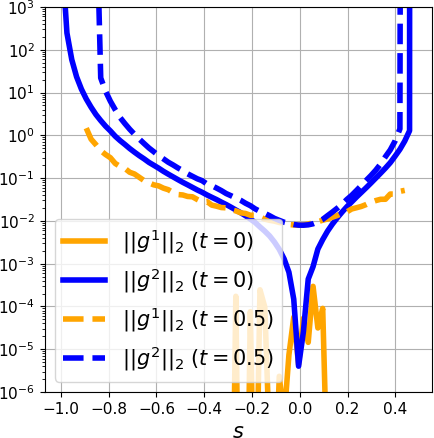}
        \hspace{4mm}
        \includegraphics[width = 0.43\textwidth]{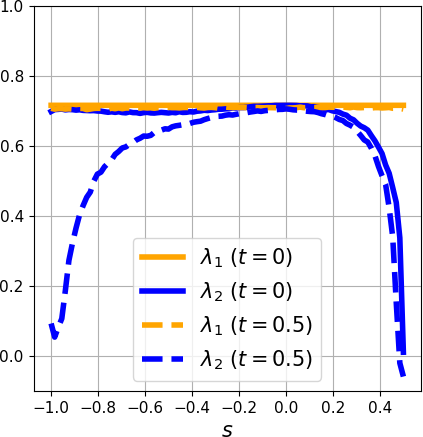}
    \end{minipage}%
    \vspace{4mm}
    \begin{minipage}{1.0\textwidth}
        \includegraphics[width = 0.45\textwidth]{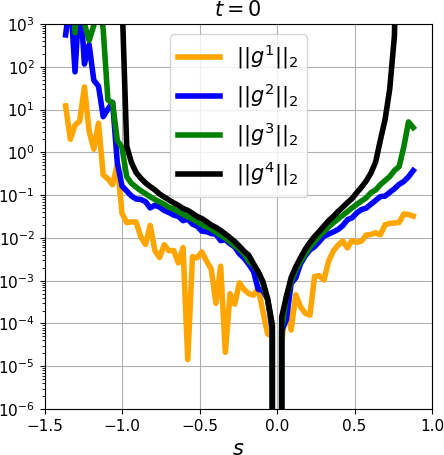}
        \hspace{4mm}
        \includegraphics[width = 0.43\textwidth]{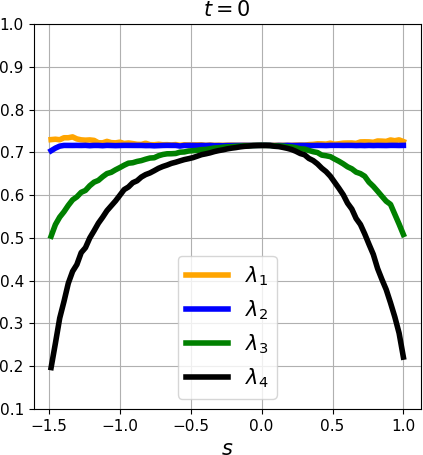}
    \end{minipage}%
    \vspace{4mm}
    \begin{minipage}{1.0\textwidth}
        \includegraphics[width = 0.45\textwidth]{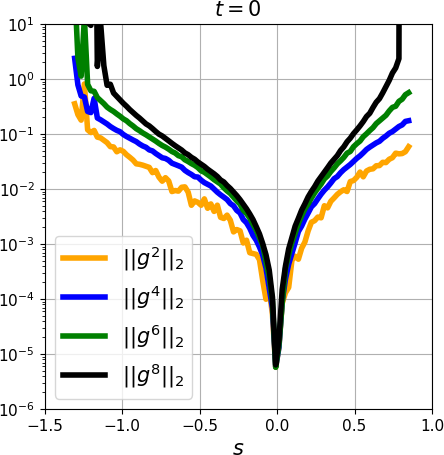}
        \hspace{4mm}
        \includegraphics[width = 0.43\textwidth]{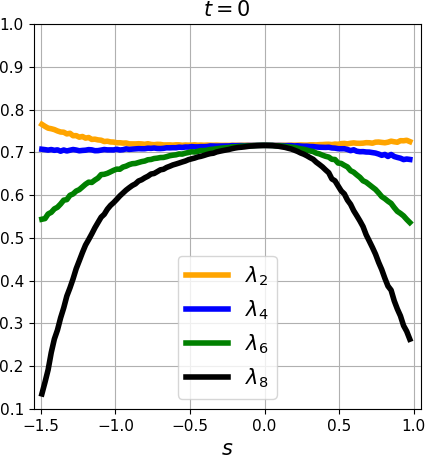}
    \end{minipage}%

    \caption{$L^2$ norm of the SRB density gradient and Lyapunov exponents of the 2D ($n = 2$; top row), 4D ($n = 4$; middle row), and 8D ($n = 8$; bottom row) variant of the coupled sawtooth map. All quantities were computed for a uniform grid of 100 values of the coupling parameter $s$. For each grid point, a trajectory of length $N = 3\cdot 10^4$ was generated.}
    \label{fig:cs-g-norm-les}
\end{figure}

In Figure \ref{fig:cs-g-norm-les}, we plot the $L^2$ norms of selected components of $g$ and corresponding Lyapunov exponents at different values of $s$ and $t$. They were computed for the 2D ($n=2$), 4D ($n=4$), and 8D ($n=8$) variants of the coupled sawtooth. In agreement with our conjecture, the norms of all components of $g$ are equal and very small in the absence of the coupling term, i.e., when $s=0$. We observe the norm ratio between $g^1$ and $g^m = g^n$ rapidly decreases as the coupling strengthens. This is also true between $g^1$ and other components corresponding to less expansive directions, as clearly indicated by the 4D and 8D examples. Figure \ref{fig:cs-g-norm-les} confirms the conjecture that the separation of Lyapunov exponents implies monotonic increase of the measure gradient norms as sorted from the most to the least expansive directions. Our results also indicate that if LEs are clustered around a single value, then the norm degradation is insignificant. Note that the converse is not necessarily true. Namely, there might be significant differences between particular components of $g$ even if LEs are clustered, which is true for the 2D sawtooth map at $s\in[-1,0]$. This usually happens when at least one of the components of $g$ is no longer integrable with respect to $\mu$ \cite{sliwiak-differentiability}. We also acknowledge the fact that square-integrability of $g$ with respect to $\mu$ is not required for the existence of linear response, as we discussed in Section \ref{sec:s3-critical}.     

In light of the specific behavior of the SRB density gradient and our main conjecture presented above, we shall numerically investigate the impact of the objective function $J$ on the statistics and their change with respect to parameters. The purpose of this experiment is to visualize long-time averages computed at different parameter values for the 2D coupled sawtooth. A fundamental question we need to raise concerns the alignment requirement. How can we say that a chosen $J$ is in fact aligned with $q^1$? Indeed, the two components of the corresponding vector $Z$ generally depend on both phase space coordinates. In the 2D setting, it is relatively straightforward to find a vector field $Z$ that satisfies that requirement. If $q^1$ is approximately parallel to $[1,1]^T$ and both components of $Z$ depend on $x^1+x^2$ only, i.e., $Z=Z(x^1+x^2)$, the corresponding $J$ is automatically aligned with $q^1$. However, if $Z^1 = Z^1(x^1+x^2)$ and $Z^2 = Z^2(x^1-x^2)$, then their respective $L^2$ norms are expected to be similar. Finally, if $Z=Z(x^1-x^2)$, $Z^2$ becomes dominant giving more weight to the second component of $g$, which is in fact the least desired scenario.  

\begin{figure}
    \centering
    \includegraphics[width = 0.54\textwidth]{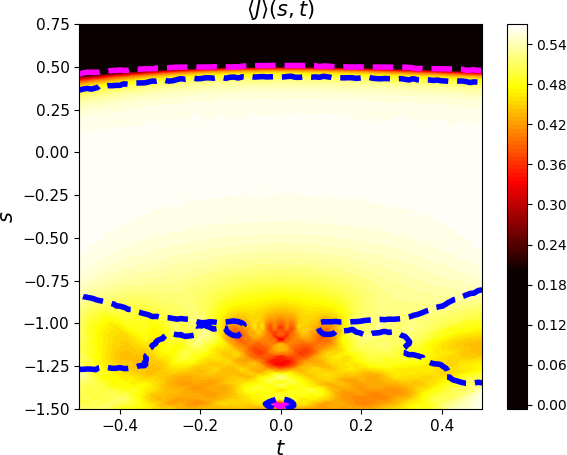}
    %\vspace{2mm}
    \includegraphics[width = 0.55\textwidth]{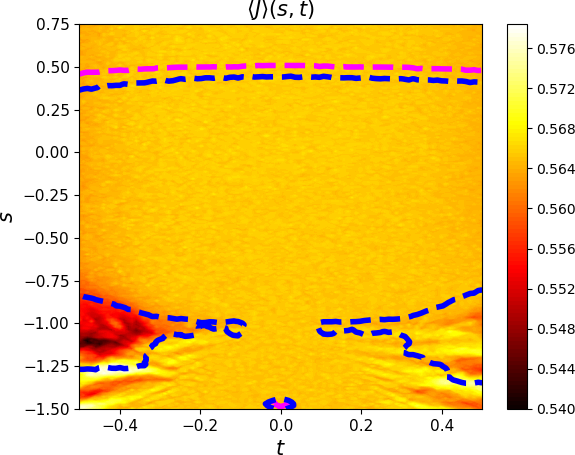}
    %\vspace{2mm}
    \includegraphics[width = 0.57\textwidth]{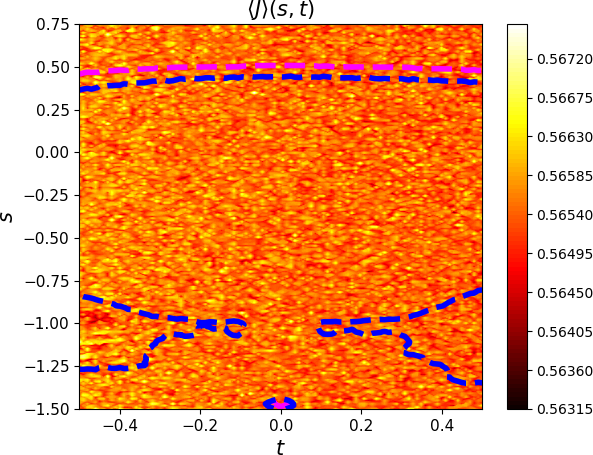}
    \caption{Long-time averages of the wave-like objective function $J=\exp(\sin(z))\,\sin(z)$, where $z = x^1-x^2$ (upper-left plot), $z = x^1$ (upper-right plot) and $z = x^1+x^2$ (lower plot). The time averages were computed for a uniform parametric grid consisting of 225 and 100 points along $s$ and $t$, respectively. For each set of parameters, a trajectory of length $N = 5\cdot10^6$ was generated. The dashed lines represent isolines corresponding to two different values of the second (i.e., smaller) LE: 0.5 (dark blue) and 0 (violet).}
    \label{fig:cs-stats}
\end{figure}

Thus, we shall consider three wave-like objective functions that depend on $x^1-x^2$, $x^1$, and $x^1+x^2$. These waves have zero gradients in the phase space directions parallel to $[1,1]^T$, $[0,1]^T$ and $[1,-1]^T$. They respectively represent functions that are weakly, moderately, and strongly aligned with the most expansive direction of the 2D hyperchaotic map. The statistics corresponding to these objective functions evaluated at a fine parametric grid are plotted in Figure \ref{fig:cs-stats}. We observe that the variation of statistics of $J = J(x^1-x^2)$ is quite large in the regions that coincide with the parametric regime of large measure change. Within this parametric subset, the value of the second LE evidently decreases and approaches the value of zero. Indeed, the largest sensitivity of the system is observed as $s$ increases from $s\approx 0.35$ to $s\approx 0.5$ for all $t\in[-0.5,0.5]$. Thus, for this parametric regime, the maximum value of $|d\langle J\rangle/ds|$ is $\mathcal{O}(1)$. In the moderate case, variations of $\langle J\rangle$ are significantly smaller compared to the previous example. However, we still observe non-negligible sensitivities of order $\mathcal{O}(10^{-1})$ if $s<-0.75$ and $|t|>0$. The third plot of Figure \ref{fig:cs-stats} shows the statistics of a function that is aligned with the most expansive direction, i.e., it depends on $x^1+x^2$. The computed long-time averages now oscillate between two values that are $\mathcal{O}(10^{-3})$ apart, across the entire parametric space. These oscillations are distributed uniformly, even around the regions of large measure gradients and distortions. In this case, $\langle J \rangle$ is approximately independent of both parameters, which implies negligible linear response.         

%Write major conclusions and implications 
The major conclusion that follows from the above analysis and numerical examples is that the unstable part of linear response might be negligible for a particular class of objective functions $J$. This is true for any system parameter with respect to which the sensitivity is computed. We observed that a scattered distribution of the positive part of the LE spectrum leads to the norm increase of consecutive components of the SRB measure gradient, represented by $g$. This usually causes significant variations of the statistics in the parameter space and, simultaneously, enables finding the optimal alignment of $J$. In this section, we demonstrated that the elimination/neutralization of the largest components of the SRB measure gradient might dramatically reduce the unstable contribution. This can be achieved by choosing a $J$ that is aligned with the most expansive direction, which is reflected by the partial integration in Eq. \ref{eqn:approx-unstable-parts}. In high-dimensional systems, we expect substantial reductions of the unstable contribution as long as $J$ is aligned with any subspace spanned by the most expansive directions. Note also that our argument applies only to systems with at least two positive LEs. If $m=1$, there is only one expansive direction, which means there are no degrees of freedom for choosing an appropriate $J$.  

How can these results and analysis be used in the context of practicable high-dimensional systems? In a standard engineering design process, the quantity of interest is a well-defined function with a concrete physical meaning, e.g., temperature, kinetic energy, drag force, that is generally not aligned with some abstract subspace of the chaotic attractor. In the following section, we argue that the specific condition imposed on $J$ is not an obstacle for a vast family of dynamical systems encountered in many fields such as climate science and turbulence theory. We show that the stable part alone can approximate the total linear response sufficiently well.  

%------------------------Higher-dimensional systems-------------------------

\section{Sensitivity analysis of higher-dimensional flows with statistical homogeneity}\label{sec:appproximating-hd}
%Introduce the concept of statistical homogeneity and explain why it is a game changer. Give examples of such systems
We presented an argument supporting the concept of small unstable contributions. This promising observation may lead to a significant simplification of the S3 algorithm for linear response. As described in Section \ref{sec:appproximating-unstable}, the major requirement for the leading unstable term $U$ to be small is a concrete alignment of the objective function $J$. In an ideal setting, the slope (variation) of $J$ in the least expansive directions should be relatively low compared to the most expansive one represented by $q^1$. This requirement seems to be very restrictive given complicated dynamical behavior of general high-dimensional chaos. In the simple example introduced in Section \ref{sec:appproximating-unstable}, the most expansive direction was predictable, thanks to which one could easily choose a suitable $J$. In this section, we will focus on a common feature of a vast group of spatially-extended chaotic systems: statistical homogeneity in space. Relying on this property, we argue that the system's dimension $n$ increases the probability of desired alignment, regardless of the physical meaning and form of $J$.   

Statistical homogeneity in the physical space implies that the long-time behavior of all system coordinates is approximately the same. For such systems, the objective function is usually defined in terms of the spatial average of a physical quantity. For 1D-in-space continuous systems bounded by $a\in\mathbb{R}$ and $b\in\mathbb{R}$, $b > a$ , for example, $J$ is usually expressed as follows,
\begin{equation}
    \label{eqn:approx-intro1}
    J = \frac{1}{b-a}\int_a^b \tilde{J}(x)\,dx \approx \frac{1}{n}\sum_{i=1}^n \tilde{J}(x^i) := \frac{1}{n}\sum_{i=1}^n \tilde{J}^i
\end{equation}
where $\tilde{J}:\mathbb{R}\to\mathbb{R}$ is a function with a concrete physical meaning. In case of the Navier-Stokes model, $\tilde{J}$ is linear if the velocity is the quantity of interest. For energy-like quantities, such as the kinetic energy, $\tilde{J}$ could be a quadratic function. Note that if the property of statistical homogeneity holds, then
\begin{equation*}
    \langle J\rangle = \langle \tilde{J}^1\rangle = \langle \tilde{J}^2\rangle = ... = \langle \tilde{J}^n\rangle,
\end{equation*}
where $\langle \cdot \rangle$ denotes the long-time average. This implies that for any time-dependent weight vector $w(t)\in\mathbb{W}$, 
where 
\begin{equation*}
    \mathbb{W} = \bigg\{ w\in\mathbb{R}^n\,\big|\,\sum_{i=1}^n w^i(t) = 1\;\forall t\geq 0 \bigg\},
\end{equation*}
the following is true
\begin{equation}
    \label{eqn:approx-intro2}
    \langle J_w \rangle := \langle \sum_{i=1}^n  w^i\,\tilde{J}^i\rangle = \sum_{i=1}^n  \langle w^i\,\tilde{J}^i\rangle \stackrel{\mathrm{indep.}}{=} \langle \tilde{J}^1\rangle \langle \sum_{i=1}^n w^i\rangle = \langle J \rangle.
\end{equation}
Eq. \ref{eqn:approx-intro2} assumes $\tilde{J}^i$ and its corresponding weight are statistically independent. Therefore, the original objective function $J$ can be replaced by any member from the class of spatially weighted functions without affecting the long-time behavior. This critical observation implies that for any smooth $J$, the feasible space of $J_w$ increases with the system's dimension $n$. It means that for a large $n$, there might be a lot of candidates well-aligned with $q^1$. Note that $w$ should primarily depend on $q^1$, i.e., an inherent topological property of the tangent space, which justifies the assumption of statistical independence of $w$ and a single phase space coordinate and, consequently, independence of $\tilde{J}^i$ and $w$ in the limit $n\to\infty$.  

We highlight yet another common property of larger physical systems. As reported by several publications (see \cite{pazo-clvs} and references therein), one can distinguish spatially localized structures of the expansive part of the Lyapunov basis. For example, in a 3D turbulent flow past a cylinder studied in \cite{ni-jfm}, the most expansive directions tend to be localized in the areas of primary instability. These include the boundary layers and near weak regions. In far wake regions and in the free steam, $q^1$ was reported to be inactive, i.e., approximately zero. Moving away from the regions of primary instability, less expansive and contracting Lyapunov vectors tend to be dominant. However, as pointed out in \cite{pazo-clvs}, in homogeneous systems with periodic boundary conditions, the clustered activity regions of $q^1$ may move across the entire physical domain. In their analysis of Rayleigh-B\'enard convection \cite{xu-clvs}, the authors notice that, for the most expansive tangent vectors, the energy spectral density is concentrated around a specific wave number, which turns out to be approximately the same as the one of the flow field (primal solution). The same work demonstrates that the energy spectrum density gradually becomes uniform as the Lyapunov vector index increases. Based on the rich numerical evidence, we expect that any time instance of $q^1$ is expected to involve local activity patterns that are restricted to a sub-region or wobble around the entire domain.

Given these specific properties of higher-dimensional chaos, the problem of alignment of $J$ and $q^1$ could be easily circumvented. Notice that we have freedom in choosing time-dependent weights, which can potentially favor only those coordinates that correspond to the regions of ``activity" of $q^1$. As these ``activity" clusters move around in time, the corresponding weights can be adjusted accordingly keeping the remaining components of $w$ close to zero. If $\tilde{J}^i = x^i$, then the optimal choice choice of weights is strictly determined by the components of $q^1$. For higher-order polynomial objective functions, the relative values of state components would also affect the corresponding weights. Their individual contributions, however, is negligible if $n$ is large. A high density of spatial coordinates facilitates search of the optimal set of weights favoring the active components of $J$ in the right proportion, regardless of the form of $\tilde{J}^i$. For a dynamical system with arbitrary statistical behavior and complex tangent topology, it is generally difficult to analytically estimate how large $n$ should be to ensure the satisfactory alignment of $J_w$ leading to the neutralization of the unstable term. Therefore, in this section, we resort to numerical studies of systems with statistical homogeneity to guarantee that Eq. \ref{eqn:approx-intro2} holds.

Before we discuss the numerical results, we first focus on algorithmic consequences of neglecting the effect of the SRB measure change. Indeed, a complete omission of the unstable part in the computation of linear response dramatically simplifies the space-split algorithm. That term, obtained through partial integration, requires computing the SRB density gradient and derivatives of projections of tangent solutions onto the unstable-center subspace. These two ingredients require solving $\mathcal{O}(m^2)$ second-order tangent equations, which is by far the most expensive section of Algorithm \ref{alg:alg1}. Assuming $n$ is large, further simplifications can be introduced. Note that the neutral contribution involves an infinite series of $k$-time correlations of $c^0$ and $DJ\cdot f$ with the leading term
\begin{equation}
    \label{eqn:approx-neutral}
    N = \int_M c^0\,DJ\cdot f\,d\mu := \int_M (c^0\,|f|) \,DJ\cdot q_f\,d\mu,
\end{equation}
where $c^0$ is the projection of a center-stable component of the tangent solution onto the center subspace normalized by the length of $f$, as derived in Eq. \ref{eqn:s3-constraint1}. Notice that the form of $N$ is in fact identical as its unstable counterpart in its original form. Therefore, if our conjecture of small unstable contributions applies, then $N$ is also small and can be neglected in the linear response algorithm. Indeed, assuming $\lambda_{m} \approx \lambda_{m+1} = 0 \approx\lambda_{m+2}$, the statistical behavior of $q_{f}$, $q^m$, $q^{m+2}$ and, consequently, the $L^2$ norms of $DJ\cdot q_{f}$, $DJ\cdot q^{m}$, $DJ\cdot q^{m+2}$ are expected to be similar. Recall also that the projection coefficients $c^{f}$, $c^m$ and $c^{m+2}$ represent dot products of a component of $v$ and their corresponding tangent vectors. The direction of parametric deformation is generally independent of Lyapunov vectors. Based on this analysis, we conclude that if our conjecture of a small $U$ holds, then the computation of $N$ could also be neglected.      

Exclusion of both unstable and neutral terms from the full S3 algorithm leaves us with the stable term alone. The remaining part requires computing the regularized tangent solution through step-by-step orthogonal projection of the unstable-neutral component. Since $f$ is generally not orthogonal to the column space of $Q$, the original stabilizing procedure involves an assembly and inversion of the Schur complement $S$. We have directly used $f$ because it is always given at no cost and it allows for a straightforward derivation of a computable formula for the neutral part of the linear response. However, since we completely neglect that part as well, the process of regularizing the tangent solution can be simplified even further. Instead of using $f$ and then orthogonalizing the $(Q,f)$ tuple, we can solve one more first-order tangent equation and perform QR factorization of the extended tangent solution matrix. Thanks to this modification, we recursively generate the orthogonal basis of the unstable-center subspace and compute projections of $v$ onto that basis, which is equivalent to the original algorithm. This can be achieved by executing Lines 9-10 of Algorithm \ref{alg:alg1} by changing $m$ to $m_{ext}$, where $m_{ext}$ should ideally be equal to $m+1$. In practice, however, setting $m_{ext} = m + 1$ may lead to instabilities due to the potentially non-hyperbolic behavior of the system. Moreover, if $n$ is large, we rarely know the exact value of $m$. If our aforementioned conjecture of a small $N$ is valid for large systems, then we can project out a few additional components of the tangent space from $v$. Therefore, as long as $m_{ext}$ is close to $m + 1$, the penalty of these extra projections, in the context of sensitivity approximation, is expected to decrease as $n\to\infty$. The only practical consequence is that a few extra tangent equations will have to be solved, which barely influences the overall cost of the reduced algorithm assuming $m_{ext} - m\ll m$. Algorithm \ref{alg:alg_reduced} summarizes all steps required to approximate the sensitivity. This procedure was obtained by eliminating the unstable and neutral contributions from the full S3 algorithm. By-products of the S3 algorithm are Lyapunov exponents, included in the $le$ array, which we compute to supplement our discussion. Benettin in \cite{benettin-le} originally proposed this approach for approximating LEs.

\begin{algorithm}\label{alg:alg_reduced}
\SetAlgoLined
\SetKwInOut{Input}{Input}
\SetKwInOut{Output}{Output}
\Input{$N$, $K$, $T$, $n$, $m_{ext}$}
\Output{$d\langle J\rangle/ds \approx s/N$, largest $m_{ext}$ LEs:= $le/N$}

Randomly generate: $x_0$, $v_0$, $Q_0$ such that $\mathrm{size}(x_0) = \mathrm{size}(v_0) = (n,1)$, $\mathrm{size}(Q_0) = (n,m_{ext})$\;
Set $s = 0$ and $le = \mathrm{zeros}(m_{ext})$\;

\For(\tcp*[h]{main time loop}){$k = 0,...,N-1$}{

\If{$k\geq T$}
{$s := s + DJ_k\cdot v_k$\;
$le := le + \mathrm{diag}(\log(\mathrm{abs}(R_k)))$\;
}
$P_{k+1} = D\varphi_k\,Q_k$\;
QR-factorize $P_{k+1}$: $Q_{k+1}\,R_{k+1} = P_{k+1}$\; 
$r_{k+1} = D\varphi_k\,v_k + \chi_{k+1}$\;
$c_{k+1} = Q_{k+1}^T\,r_{k+1}$\;
$v_{k+1} = r_{k+1} - Q_{k+1}\,c_{k+1}$\;
Advance the iteration: $x_{k+1} = \varphi(x_k)$\;
}

\caption{Reduced space-split sensitivity algorithm for higher-dimensional chaotic flows}
\end{algorithm}

\subsection{Lorenz 96}\label{sec:appproximating-lorenz96}
\begin{figure}
    \centering
    \includegraphics[width = 0.35\textwidth]{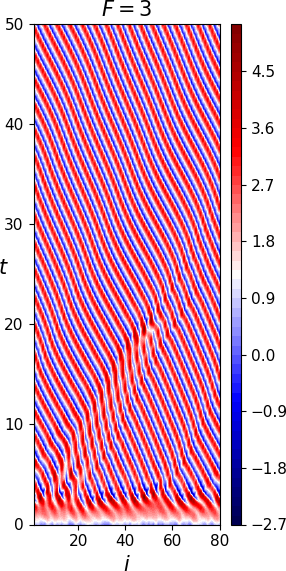}
    \hspace{1cm}
    \includegraphics[width = 0.36\textwidth]{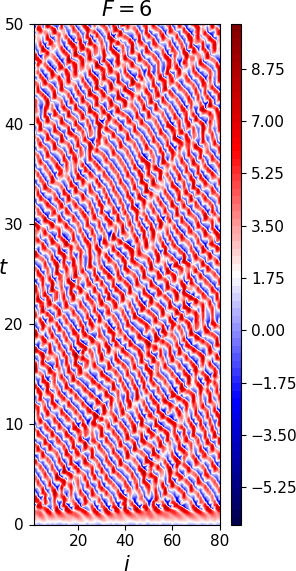}
    \vspace{3mm}
    \includegraphics[width = 0.35\textwidth]{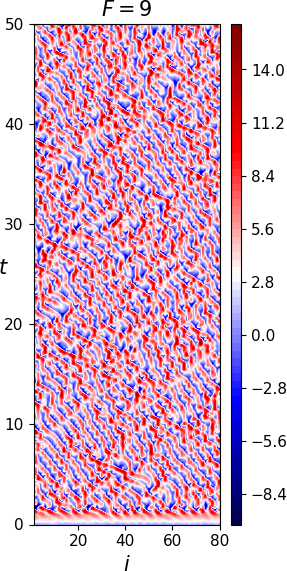}

    \caption{Solutions to the Lorenz 96 system (Eq. \ref{eqn:approx-lorenz96}) for $n=80$ stacked horizontally.}
    \label{fig:lorenz96-solutions}
\end{figure}

In light of the above conclusions, we shall consider the Lorenz 96 model, which was proposed by E. Lorenz in \cite{lorenz-96} to study spatio-temporal dynamics of the atmosphere. Mathematically, this is an $n$-dimensional chaotic flow defined as follows,
\begin{equation}
    \label{eqn:approx-lorenz96}
    \begin{split}
    \frac{dx^i}{dt} &= (x^{i+1} - x^{i-2})\,x^{i-1} - x^{i} + F,\;\;\;i = 1,...,n,\\
    x^{i+n} &= x^{i},
    \end{split}
\end{equation}
where the superscript indicates the component index, in compliance with our notation convention. Each degree of freedom $x^{i}$ represents a value of a physical quantity, e.g., temperature or pressure, on a uniformly discretized parallel of the Earth. Analogously to semi-discretized PDEs describing advection, this system involves spatially coupled variables with a quadratic nonlinearity. Eq. \ref{eqn:approx-lorenz96} involves two constant parameters: the number of sectors $n\geq 4$, each corresponding to a different meridian of the Earth, and imposed forcing $F\in\mathbb{R}^{+}$. If $F<8/9$, then the solution quickly decays to the constant value of $F$, i.e., $x^{i}=F$, $i=1,...,n$ for all $t>t^{*}\approx 0$ \cite{karimi-lorenz96}. We solve Eq. \ref{eqn:approx-lorenz96} using the explicit fourth-order Runge-Kutta with $\Delta t = 0.005$. This ODE solver will be used throughout this section, unless stated otherwise. In Figure \ref{fig:lorenz96-solutions}, we plot the solutions for $n = 80$ and three different values of $F$. For $F=3$, the periodic dynamics involves waves travelling to the west (i.e., in the direction of decreasing sector index $i$). The distortion that appears at the beginning of the simulation quickly decays leading to a predictable behavior. While some regularity is still maintained at $F=6$, the alignment of waves seems random which implies that some unstable modes might be activated. If we further increase $F$ to the value of 9, the spatio-temporal structure of the solution clearly reflects chaotic behavior without any distinguishable patterns.  

\begin{figure}
    \centering
    \includegraphics[width = 0.44\textwidth]{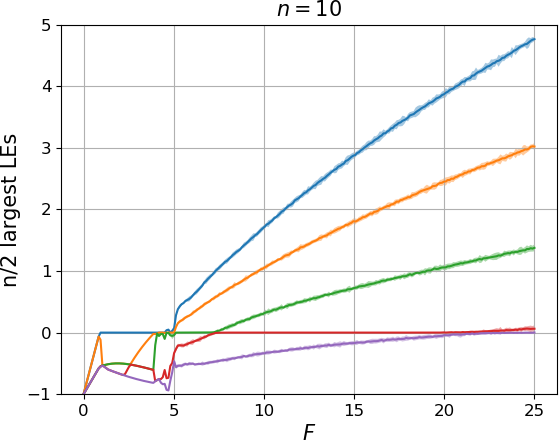}
    \hspace{2mm}
    \includegraphics[width = 0.44\textwidth]{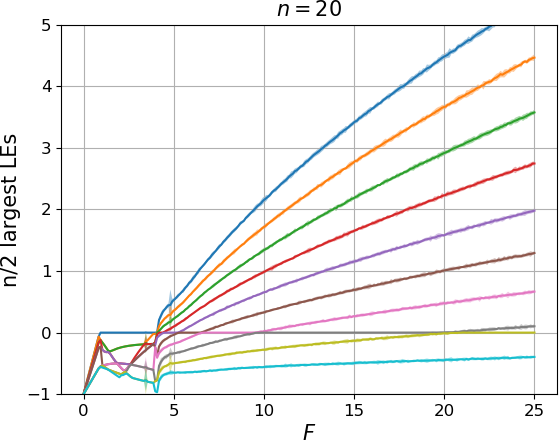}   
    
    \vspace{2mm}
    \includegraphics[width = 0.44\textwidth]{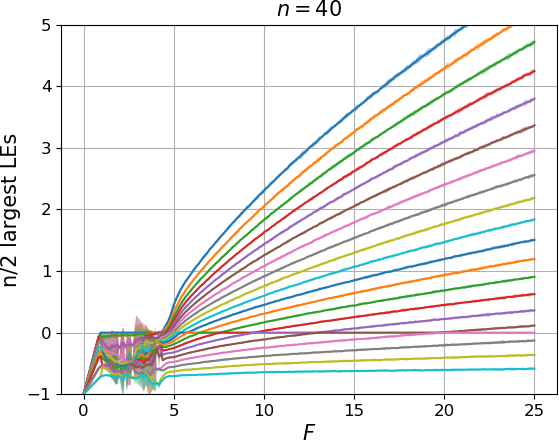}
    \hspace{2mm}
    \includegraphics[width = 0.44\textwidth]{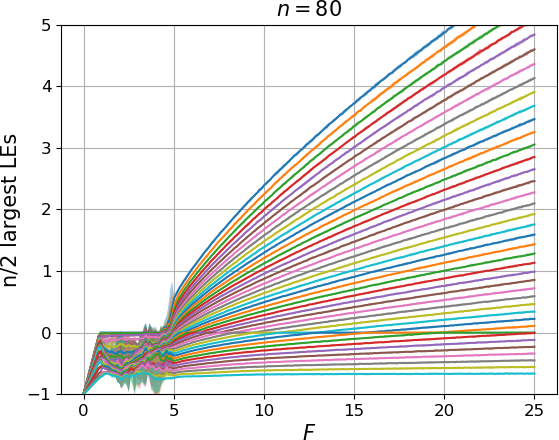}
    \caption{Larger half of the Lyapunov spectrum of Eq. \ref{eqn:approx-lorenz96}. LEs were computed at 240 distinct values of $F$ distributed uniformly between 0 and 25. For each value of $F$, we run 10 independent simulations for 5000 time units. The barely visible shaded area represents the 2-sigma range (95\% confidence) of the 10-element data set at each value of $F$.}
    \label{fig:lorenz96-les}
\end{figure}

To obtain more insights into the dynamics of the Lorenz 96 model, we analyze its Lyapunov spectrum for the most common values of system's parameters \cite{kekem-lorenz96}. In Figure \ref{fig:lorenz96-les}, we illustrate a half of the Lyapunov spectrum for $F\in[0,25]$ at $n = 10, 20, 40, 80$. For any $n$ and $F < 0.9$, all LEs are negative, which means that, for any random initial condition, the solution exponentially decays to a constant value. Within the interval $F\in[0.9,4.5]$, the dynamics is no longer stationary, but still non-chaotic, because $\lambda_{1} = 0$. We observe the presence of at least one positive LE if $F > 4.5$. In the chaotic regime, the dimension of the expansive manifold gradually increases with $F$ to about $m = n/2$ at $F=25$. Notice also that the higher $F$, the smaller the angle between the lines representing $\lambda_{i}(F)$, $i=1,2,...$ and the x-axis. Indeed, the authors of \cite{karimi-lorenz96} computed a curve fit for $\lambda_1^{-1}(F)$ at $N=35$, whose close-form formula is the following: $\lambda_1^{-1}(F) = 0.158 + 123.8\,F^{-2.6}$. Consequently, given the self-similar behavior of the plotted spectrum, all LEs seemingly converge to fixed values as the forcing $F$ increases.      

\begin{figure}
    \centering
    \includegraphics[width = 0.7\textwidth]{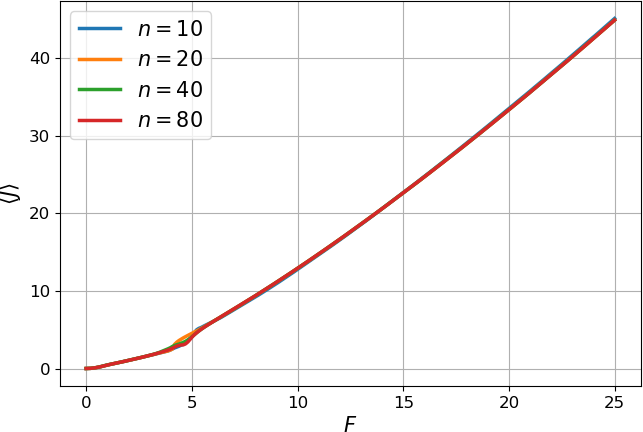}
    \caption{Long-time means of spatially-averaged kinetic energies of the Lorenz 96 system. The statistics were computed on a uniform grid of 240 values of $F\in[0,25]$. For each value of $F$, the objective function was time-averaged over $5\cdot 10^6$ time units.}
    \label{fig:lorenz96-jav}
\end{figure}

We shall consider the spatially-averaged kinetic energy of the system as the objective function $J$, which can be expressed using Eq. \ref{eqn:approx-intro1} with $\tilde{J}^i = (x^i)^2$. The long-time averages $\langle J\rangle$ for $F\in[0,25]$ at $n = 10,20,40,80$ are plotted in Figure \ref{fig:lorenz96-jav}. We observe that all four curves $\langle J\rangle (F)$ collapse into a single curve due to spatial averaging. The only misalignment occurs at the non-chaotic/chaotic transition region close $F = 5$. Thus, in the extensive chaos regime of Lorenz 96, the spatially-averaged statistics is generally independent on $n$, which was previously observed in \cite{karimi-lorenz96}. We shall restrict our attention to that regime, i.e., when $F\geq 5$, and compute sensitivities with respect to $F$ using our reduced S3 algorithm. The slope of $\langle J\rangle (F)$ seems to be constant and is approximately $2$ for $F\in [5,25]$. We will use a higher-order interpolation of the statistics curve and differentiate it using the central finite-difference scheme. This estimate will serve as a reference solution to evaluate the performance of Algorithm \ref{alg:alg_reduced}.

\begin{figure}
    \centering
    \includegraphics[width = 0.55\textwidth]{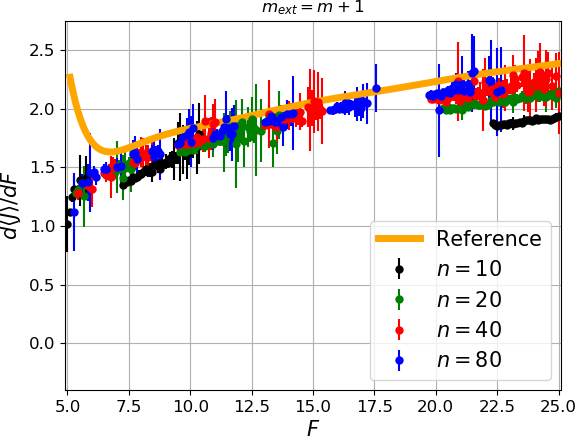}
    
    \vspace{1.5mm}
    \includegraphics[width = 0.55\textwidth]{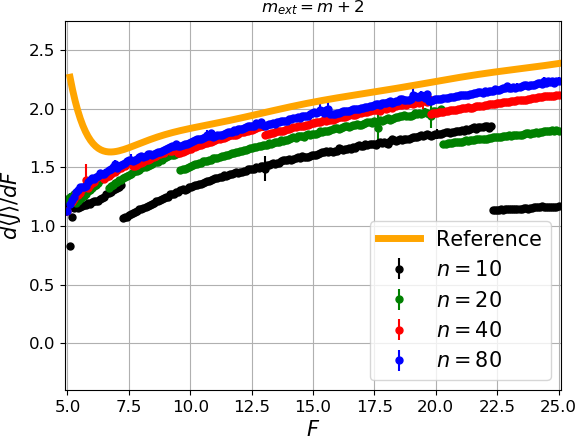}
    
    \vspace{1.5mm}
    \hspace{-1mm}\includegraphics[width = 0.55\textwidth]{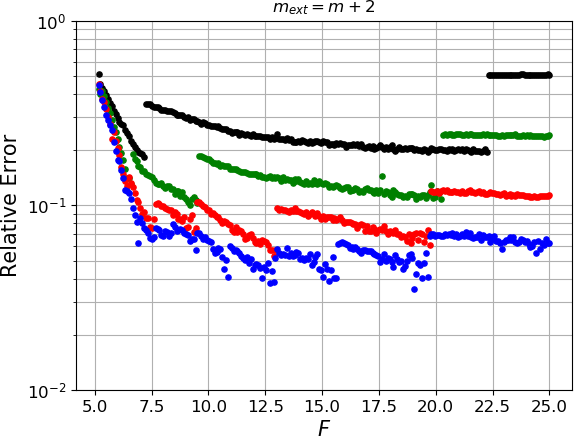}
    \caption{Linear response approximations of the Lorenz 96 model with respect to $F$ computed using Algorithm \ref{alg:alg_reduced}. The top plot illustrates computed sensitivities for $m_{ext} = m + 1$, the middle plot for $m_{ext} = m + 2$, while the bottom plot depicts the mean relative error of the $m_{ext} = m + 2$ case computed with respect to the reference finite difference solution (respective colors indicate $n$). Sensitivities were computed on a uniform 240-point grid between $F = 5$ and $F = 25$. For each value of $F$, we run 10 independent ergodic-averaging simulations over $N\Delta t = 5000$ time units. Vertical lines represent sigma intervals, while the bullets indicate the corresponding averages. Lack of a bullet (in the upper plot) means the standard deviation is larger than 1. The solid orange line is a finite difference approximation of the 11-th degree polynomial fit of $\langle J \rangle$.}
    \label{fig:lorenz96-sens}
\end{figure}

Figure \ref{fig:lorenz96-sens} illustrates approximations of linear response obtained with Algorithm \ref{alg:alg_reduced}. In particular, we used our reduced algorithm to approximate $d\langle J\rangle/dF$ for $F\in[5,25]$. For $m_{ext} = m+1$, the algorithm generates satisfactory approximation for $F\geq 6$. However, the standard deviation is quite large and it very often exceeds the value of one across the entire parametric domain. These statistical fluctuations are eliminated by increasing $m_{ext}$. Indeed, the $m_{ext} = m+2$ case has dramatically smaller sigmas everywhere. This result indicates that if $m_{ext}$ is too small, the regularized tangent solution may still have some rapidly growing components in some parts of the attractor leading to large variances. The smooth behavior of linear response in the $m_{ext} = m+2$ case suggests that these fluctuations are not caused by the ergodic-averaging error. As expected, there is always an extra penalty for increasing $m_{ext}$. However, the higher $n$, the smaller price must be paid for extra stabilizing projections. This observation is consistent with our conjecture suggesting that the relative contribution of a single component of $v$ decreases as $n$ gets larger.

Figure \ref{fig:lorenz96-sens} reveals two other critical features of the reduced algorithm. First, if $n$ is sufficiently large, then the obtained sensitivity approximation might be very accurate, i.e., the relative error within a few percent. This result confirms our major conjecture of negligible unstable (and neutral) contributions to linear response. For Lorenz 96, the impact of the SRB measure change is apparently insignificant. The only exception is the region around $F=5$. Indeed, the error is large in this parametric regime, regardless of the value of $m_{ext}$ and system's dimension $n$. Although the property of spatial homogeneity is unaffected and some unstable modes are still active, we observe the sensitivity approximation clearly deviates from the reference solution. Note this parametric region coincides with rapid value decrease of positive LEs. Many of them are still positive but there are clustered. Our discussion in Section \ref{sec:appproximating-unstable} suggests that in this case there might be no gain due to the alignment of $J$ and $q^1$. All components of $g$ are expected to have similar distributions across the phase space. Therefore, even if $J$ and $q^1$ are aligned, the unstable contribution could be significant in this case.  

\begin{figure}
    \centering
    \includegraphics[width = 0.47\textwidth]{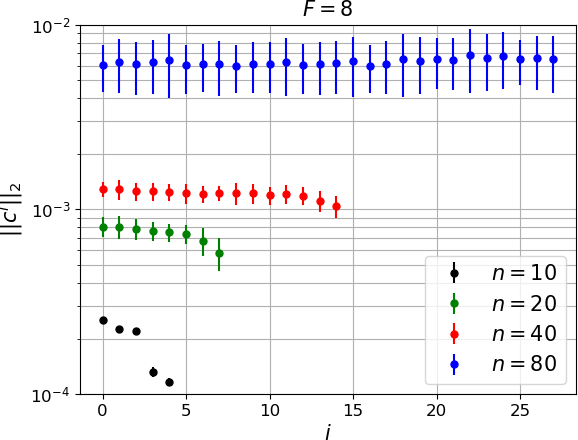}
    \hspace{2mm}
    \includegraphics[width = 0.47\textwidth]{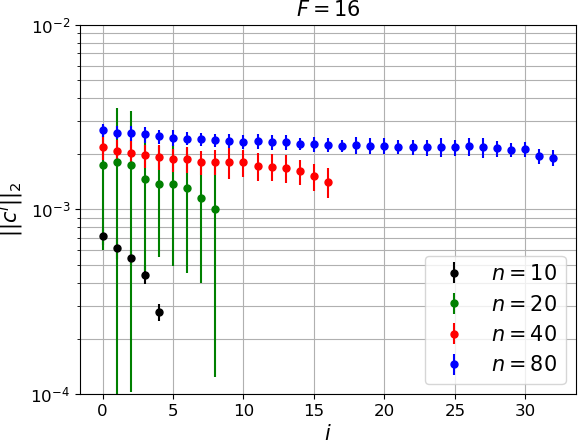}
    
    \caption{$L^2$ norms of $c^i$, $i=1,...,m_{ext}=m+2$, which were computed as by-products of Algorithm \ref{alg:alg_reduced}. All simulation parameters are the same as those reported in the caption of Figure \ref{fig:lorenz96-sens}.}
    \label{fig:lorenz96-c0}
\end{figure}

For completeness, in Figure \ref{fig:lorenz96-c0}, we also plot the $L^2$ norms of the projection scalars $c^i$, $i=1,...,m_{ext}=m+2$. This result confirms that all scalars contribute almost equally to linear response suggesting that their relative significance is degraded as $n$ increases. These results also indicate that if $n$ is small, the scalars corresponding to the lowest indices tend to be statistically larger compared to their counterparts. In other words, the Lorenz 96 system with few degrees of freedom tends favors the contributions of $\|c^i\|_2$ corresponding to the most expansive directions (small $i$).  

\subsection{Kuramoto-Sivashinsky}\label{sec:appproximating-ks}

Finally, we shall consider the Kuramoto-Sivashinsky (KS) equation, one of the simplest partial differential equations modeling chaos. Similarly to Lorenz 96, KS is a spatio-temporal description of complex dynamics driven by instabilities far from an equilibrium. This equation was proposed decades ago to model wave propagation in reaction-diffusion systems \cite{kuramoto-ks} and hydrodynamic instabilities of laminar flames \cite{sivashinsky-ks}. A number of other applications of the KS equation can be found in the literature. In this work, we analyze a modified version of KS, which includes an extra advection term proportional to a constant scalar $c\in\mathbb{R}$. The modified equation, which was previously studied in \cite{blonigan-ks}, has the following form,
\begin{equation}
    \label{eqn:ks}
    \begin{split}
    &\frac{\partial u}{\partial t} = -(u+c)\,\frac{\partial u}{\partial x} - \frac{\partial^2 u}{\partial x^2} - \frac{\partial^4 u}{\partial x^4},\\
    &u(0,t) = u(L,t) = 0,\\
    &\frac{\partial u}{\partial x}(0,t) = \frac{\partial u}{\partial x}(L,t) = 0,
    \end{split}
\end{equation}
where $x\in[0,L]$, $L=128$, $t\geq 0$, $u(x,t)\in\mathbb{R}$. We discretize this system in space using the finite difference method with second-order accuracy. The grid is uniform and involves $513$ nodes, which gives us a constant spacing $\Delta x = 128/(513-1) = 0.25$. A combination of center and one-sided schemes is applied to approximate all spatial derivatives as suggested in \cite{blonigan-ks}. The number of ODEs, i.e., the system's dimension, is reduced to $n = 511$ by incorporating all boundary conditions using the ghost node technique. While this is a stiff system, we apply the fully-explicit fourth-order Runge-Kutta scheme with a small time step $\Delta t = 0.0006$. In \ref{sec:derivatives}, we discuss how the linear response algorithm could be integrated with implicit schemes. 

\begin{figure}
    \centering
    \includegraphics[width = 1\textwidth]{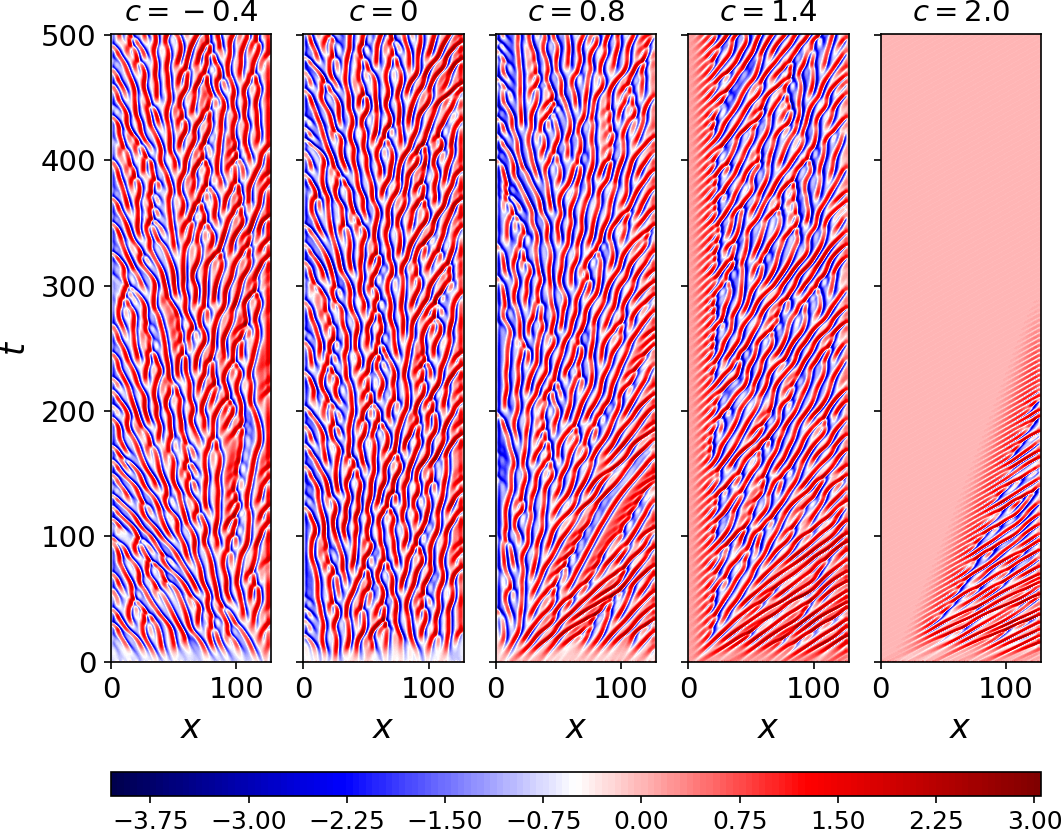}
    \caption{Solutions to the KS equation (Eq. \ref{eqn:ks}) for different advection intensities.}
    \label{fig:ks_solutions}
\end{figure}

Figure \ref{fig:ks_solutions} illustrates solutions to the KS equation, $u(x,t)$, for different values of $c$. In the spatio-temporal space, $u(x,t)$ involves a collection of irregular branches that switch between positive and negative values. The sign of $c$ determines the inclination of these branches. If $c$ is positive, they tend to move in the positive direction of $x$ and vice versa. By increasing the magnitude of $c$, the advection term starts to dominate pushing the lightly turbulent region out of the domain. Indeed, for $c=2$, we observe that $u(x,t)$ quickly becomes steady suggesting that all unstable modes are killed due to the strong advection. Regardless of the value of $c$, one can distinguish a transitional period at the beginning of each simulation during which the spatio-temporal branches develop their shapes. At $c=1.4$, the spatial sub-region $x<20$ is dominated by the convection, which results in an almost stable behavior of $u(x,t)$ in that part of the domain. This leads to violation of statistical homogeneity along $x$.      

\begin{figure}
    \centering
    \includegraphics[width = 0.89\textwidth]{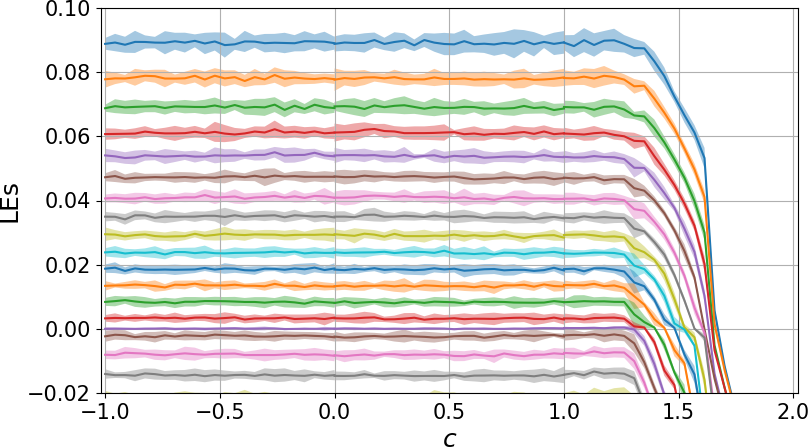}
    \caption{18 largest Lyapunov exponents of the KS equation. The spectrum was computed at the uniform grid between $c=-1$ and $c=2$. For each value of $c$, 10 independent simulations were run. The sought-after quantities were obtained through ergodic-averaging over $12,000$ time units per simulation. The solid lines represent the mean values obtained in 10 simulations, while the shaded area represents the 2-sigma range.}
    \label{fig:ks_les}
\end{figure}

Figure \ref{fig:ks_les} depicts the $18$ largest Lyapunov exponents of the KS equations for $c\in[-1,2]$. The LE spectrum is independent of $c$ as long as $-1\leq c \leq 1.3$. At $1.3 \leq c \leq 1.7$, we observe a rapid decrease of all positive LEs. This coincides with the increasing strength of the advection term. Intuitively, the dominating advection term gradually kills the unstable modes, which consequently leads to a more predictable behavior of $u(x,t)$. The KS system is clearly non-chaotic if $c>1.7$, which is reflected by the stable behavior of $u(x,t)$ at $c=2$ illustrated in Figure \ref{fig:ks_solutions}.

We also acknowledge similarities in the behavior of LE spectra corresponding to the Lorenz 96 and KS system. In the former, we observed an analogous collapse of the values of positive LEs around the stability-to-turbulence transition close to $F=5$. Another analogy is the parametric independence of the LE spectrum at large values of $F$. Note, however, that the ratio $m/n$ may reach the value of $1/2$ in case of Lorenz 96, which is significantly larger compared to this case.

\begin{figure}
    \centering
    \includegraphics[width = 0.9\textwidth]{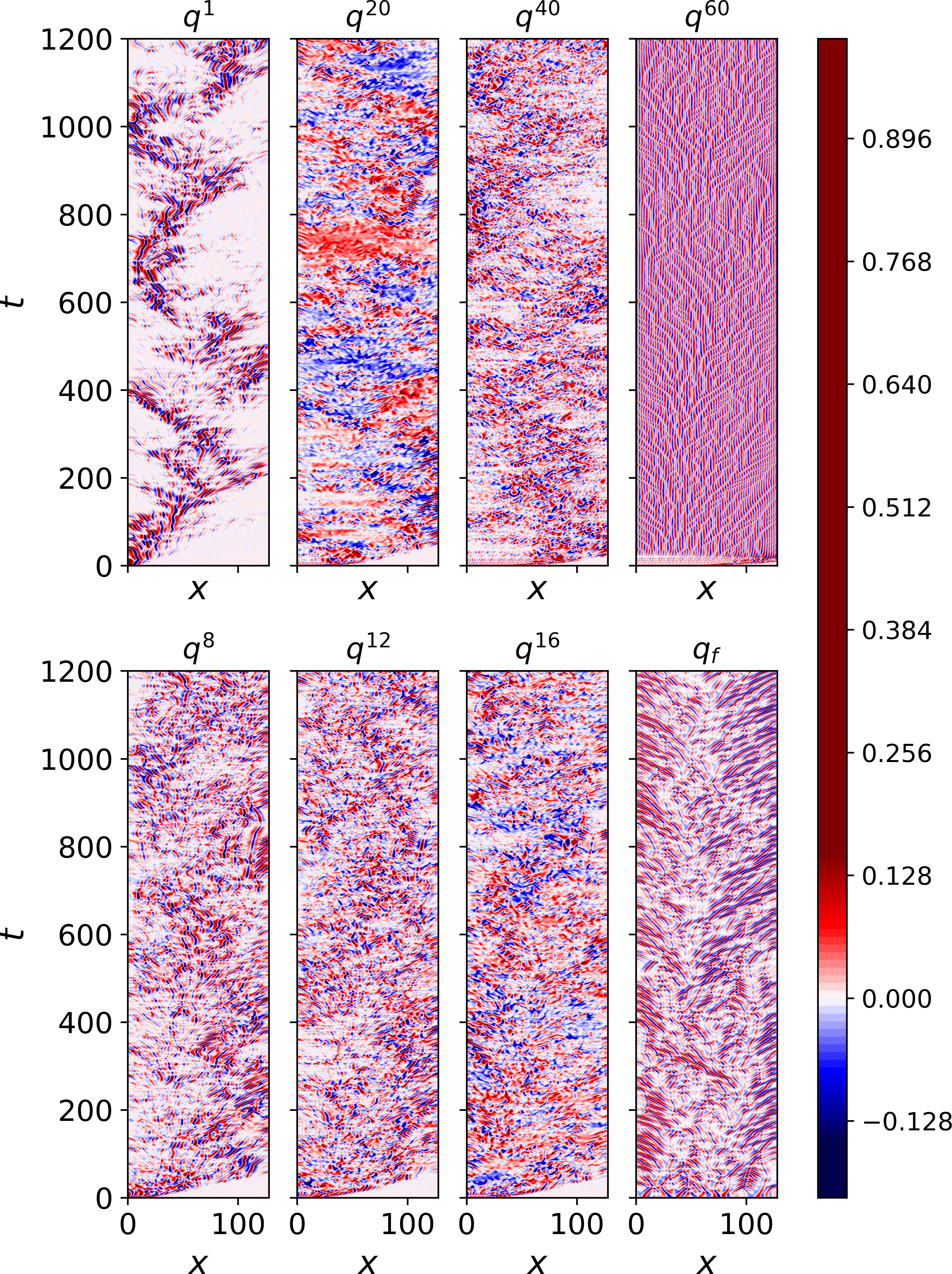}
    \caption{Orthonormal Lyapunov vectors $q^i$ of the KS system (Eq. \ref{eqn:ks}) without the extra advection term ($c=0$). The vector $q_f$ represents the normalized time derivative of $u(x,t)$. The colorbar has linearly been re-scaled between $-0.15$ and $0.15$ keeping the same color for all values from beyond this interval.}
    \label{fig:ks_laps}
\end{figure}

Selected Lyapunov vectors are plotted for $t\in[0,1200]$ in Figure \ref{fig:ks_laps}. As expected, the leading Lyapunov vector $q^1$ consists of relatively large structures with local support. The region of activity of $q^1$, which corresponds to non-small components, is limited to a thin sub-region, which moves around the entire $x$-space. It periodically bounces back and forth between the two walls. We observe that the structural behavior of $q^i$ visibly changes as $i$ increases. The $q^{40}$ vector features much finer structures with occasional inactivity regions, while $q^{60}$ seems to be periodic and highly-oscillatory in $x$, and almost constant in $t$ across the entire domain. It is quite surprising that $q^{20}$ features large compact shapes, such as the red one at $t\in[650,800]$. The tangent vectors corresponding to less expansive and mildly contracting directions are placed in the bottom row of Figure \ref{fig:ks_laps}. They consist of finer structures compared to the ones of $q^1$ and have occasional small inactivity regions throughout the entire domain. All vectors in the bottom row are visibly similar except when $t$ is small. Indeed, all Lyapunov vectors $q^i$ were obtained in an iterative procedure that is initiated at a random initial condition. We observe that this iteration requires at least 50 time units for a run-up.       

\begin{figure}
    \centering
    \includegraphics[width = 0.49\textwidth]{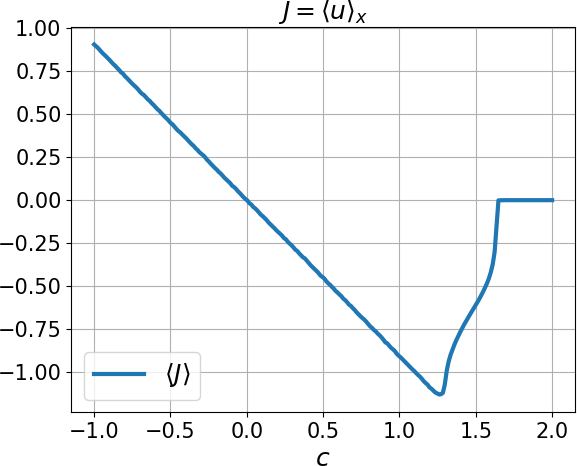}
    \hspace{2mm}
    \includegraphics[width = 0.47\textwidth]{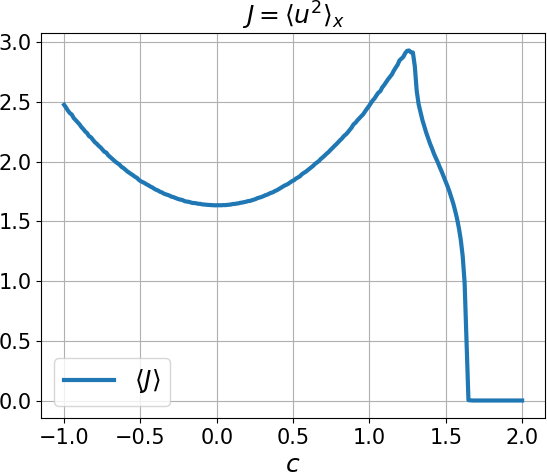} 
    
    \includegraphics[width = 0.47\textwidth]{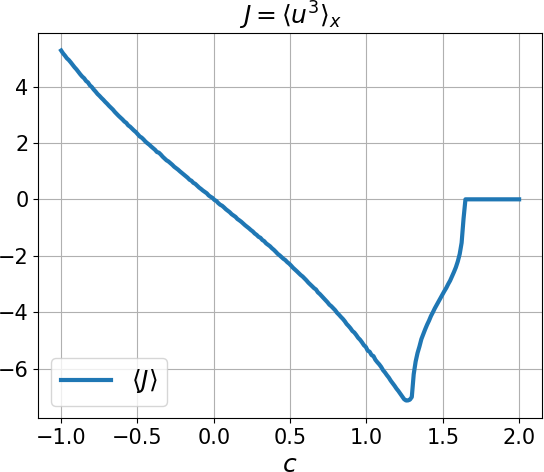} 
    
    \caption{Long-time averages $\langle J\rangle$ computed on a uniform 240-point grid of $c\in[-1,2]$. The operator $\langle\cdot\rangle_x$ indicates the spatial average. For each value of $c$, we run an ergodic-averaging simulation over 600,000 time units.}
    \label{fig:ks_stats}
\end{figure}

Given these preliminary results, we apply Algorithm \ref{alg:alg_reduced} to compute linear response with respect to the parameter $c$. This time we shall consider three different spatially-averaged objective functions: linear, quadratic and cubic, i.e., $\tilde{J}^i = u^{p}$, $p = 1,2,3$, respectively. The corresponding long-time averages are plotted in Figure \ref{fig:ks_stats} at $c\in[-1,2]$. We observe that in all of these cases the mean curve can be divided into three smooth sections connected at $c\approx 1.25$ and $c\approx 1.7$. The shape of the left part resembles a polynomial function of the same order as the objective function itself. The middle one resembles the tangent function, while the right-hand side piece is constant in all three cases. These three pieces coincide with three different behavior types of $u(x,t)$ that we observed in Figure \ref{fig:ks_solutions}: turbulent ($c\leq 1.25$), transitional ($1.25\leq c\leq 1.7$), and advection-dominated ($c\geq 1.7$) regime.

\begin{figure}
    \centering
    \includegraphics[width = 0.75\textwidth]{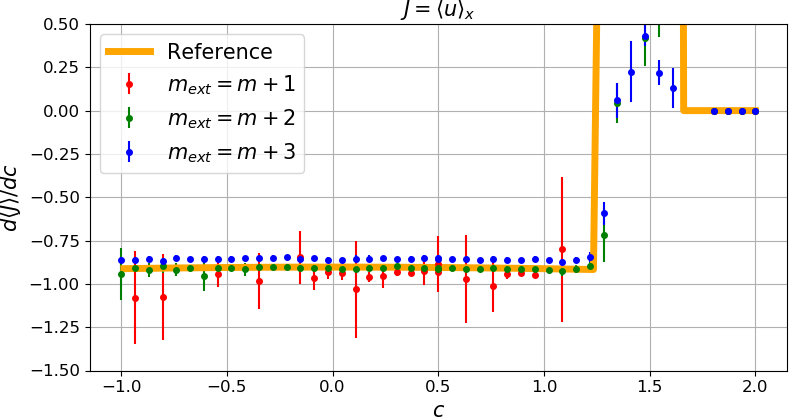}
    
    \vspace{5mm}
    \includegraphics[width = 0.75\textwidth]{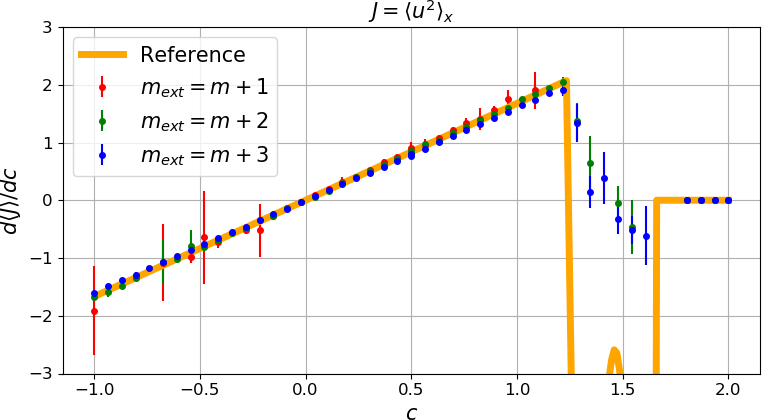} 
    
    \vspace{5mm}
    \includegraphics[width = 0.75\textwidth]{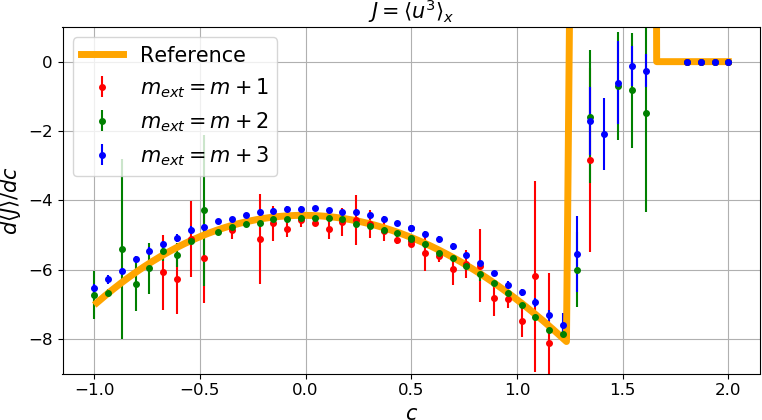} 
    
    \caption{Linear response computed for the same objective functions as those presented in Figure \ref{fig:ks_stats} using Algorithm \ref{alg:alg_reduced}. For each value of $c$, we run 10 independent simulations over 3,000 time units each. Bullets and vertical lines represent the mean and standard deviation, respectively. The results with a large standard deviation were removed from the plot. The reference line was computed through central finite-differencing of polynomial fits.}
    \label{fig:ks_sens}
\end{figure}

We apply our reduced linear response algorithm (Algorithm \ref{alg:alg_reduced}) to approximate sensitivities for these three objective functions. Analogously to the previous plots, we compare our approximation against the finite-difference reference solutions. Figure \ref{fig:ks_sens} illustrates linear response results for different values of $m_{ext}$. One can easily observe a lot of similarities between these results and the ones generated for Lorenz 96. First of all, if $m_{ext} = m+1$, then the mean solution is quite close to the reference line, but the variance is likely to be large. The variance is significantly reduced by increasing $m_{ext}$ and, in most cases, the new mean approximations are still very accurate. Indeed, the accuracy can be within the reference line width in the turbulent and stable regimes. Huge disparities occur in the transitional regime, i.e., at $c\in[1.25,1.7]$. Similarly to the Lorenz 96 case, this region corresponds to the sudden decrease of positive LEs. The approximation errors here are generally smaller compared to those computed for the Lorenz 96 system. Recall that, in Figure \ref{fig:lorenz96-sens}, we observed that the approximation error decreases as $n\to\infty$. Indeed, the dimension of the discretized KS system is an order of magnitude larger than that of Lorenz 96.  

Our numerical results presented in this section indicate that linear response of a higher-dimensional system can be accurately approximated by the reduced S3 method. That algorithm, which was simply obtained by eliminating the unstable and neutral contributions, solves a regularized tangent equation by projecting out all expansive and, sometimes, a few mildly contracting components. This process can be in fact formulated as an optimization problem in which we minimize the $L^2$ norm of the sum of the standard tangent solution and a linear combination of expansive orthogonal Lyapunov vectors. A similar concept was previously utilized in a variant of shadowing methods known as NILSS \cite{ni-nilss}, which relies on Covariant Lyapunov Vectors (CLVs). There are some algorithmic differences between the reduced S3 and NILSS, as the latter uses larger times windows to update tangent solutions, which requires solving extra linear systems to preserve continuity. This work sheds light on the reliability on relatively simple methods relying on some form of a regularized tangent equation.  

\section{Conclusions}\label{sec:conclusions}
%Summarize full S3 - pros and cons
%Alignment of J - its consequences
%Summarize main results of numerical experiments
%Major consequences of this work and future endeavors

Sensitivity analysis of chaotic dynamical flows is full of mathematical and algorithmic challenges. The linear response theory, especially Ruelle's formalism, allows us to better understand how different dynamical features of a system affect its sensitivity. In particular, we can rigorously decompose linear response into three separate ingredients: unstable, neutral, and stable. This concept has been utilized in recently developed algorithms such as the space-split sensitivity (S3). The unstable part represents the effect of the SRB measure gradient, which requires computing second derivatives of coordinate charts describing unstable manifolds and, equivalently, differentiating Lyapunov vectors in all unstable directions. The neutral and stable parts, as their names suggest, reflect the contributions of the parametric perturbation along the center (tangent to the flow) and stable manifolds, respectively. In general, any of these three terms might significantly contribute to the total response. The example of Lorenz 63 clearly indicates that neglecting the unstable or neutral term leads to large errors.

Despite their elegance, rigor and accuracy, direct linear response algorithms have certain flaws. First of all, they are expensive. The leading flop count may be proportional even to the cube of the number of positive Lyapunov exponents. In addition to that, the non-hyperbolic behavior of larger systems could cause numerical instabilities making the computation of measure gradients difficult. We observed that the most expansive components of the measure gradient tend to be significantly smaller in norm compared to the other ones. This critical observation led us to the conjecture that the unstable contribution could potentially be reduced if the effect of the larger components of the measure gradient is eliminated. To make the unstable part small, regardless of the parameter with respect to which linear response is computed, one could choose an {\it aligned} objective function $J$. We show that if $J$ is represented by the unstable divergence of a smooth vector field such that the directional derivative in the most expansive direction is dominant, then the majority of the measure gradient components could be killed. Our experiment on the hyperchaotic {\it coupled sawtooth map} confirms that the unstable part can be significantly reduced through an appropriate manipulation of $J$.

While the idea of finding an {\it aligned} $J$ may seem to be a purely theoretical concept, we argue that this result could be critical for practitioners as well. Indeed, spatially-extended high-dimensional chaotic systems with statistical homogeneity in space do allow for different representations of $J$. In particular, the objective function, which typically equals the spatial average of system coordinates (or higher-order moments), can be represented by an arbitrary linear combination of individual coordinate terms. Consequently, this gives us freedom in choosing $J$ and increases the probability of finding an {\it aligned} $J$ as the system's dimension grows. This conjecture is verified by eliminating the unstable and, consequently, the neutral part from the full S3 algorithm. Leaving the stable contribution alone, we accurately approximated sensitivities in both the Lorenz 96 and Kuramoto-Sivashinsky models.

Two primary goals were achieved in this work. First, we presented some parts of the full linear response algorithm with critical analysis of its major parts and potential applications. Second, based on our analysis, we proposed a reduced variant of S3 that has been shown to be sufficient for some higher-dimensional systems. Our results indicate that, in systems with statistical homogeneity, sensitivities could be accurately approximated by projecting out the unstable components from the tangent solution. Hence, the effect of the SRB measure change is negligible on those systems for a wide range of parameters. We showed that when the Lyapunov spectrum collapses, which typically happens when the system moves from a non-chaotic to chaotic regime, the stable term alone is not enough. Our future work shall investigate how likely this scenario is in real-world engineering applications. If this is a rare event, further developments of well-established shadowing methods would not be necessary. Otherwise, one could consider extracting some parts of the unstable contribution to correct the reduced algorithm. 

%----------------Miscellaneous-------------------

\section*{Supplementary Material}

To facilitate the reproduction of the reported results, we attach our Python code that is sufficient to generate and post-process data relevant for this paper. In the software package, the reader will find the README file with a description of the attached scripts.  

\section*{Acknowledgments}
This work was funded by U.S. Department of Energy Grant No. DE-NA-0003993. The authors also acknowledge the MIT SuperCloud and Lincoln Laboratory Supercomputing Center for providing HPC resources that have contributed to the research results reported within this paper.

\section*{Conflict of interest}
The authors declare that they have no conflict of interest.

\bibliographystyle{elsarticle-num-names}
\bibliography{references.bib}

%----------------Appendices-------------------
\appendix

\section{Full space-split algorithm -- description, pseudocode and complexity analysis}\label{sec:algorithm}

The purpose of this section is to extend the discrete version of S3 \cite{sliwiak-s3} to continuous chaos and present the structure of the full linear response algorithm. We rely on the three-term splitting defined by Eq. \ref{eqn:s3-splitting}. The major difference difference between the discrete and continuous variants of S3 is that, in the latter, we additionally project out the neutral component from the regularized tangent solution $v$. The computation of the stable part involves solving a linear system for $c^i$, $i=0,1,...,m$, because the vector tangent to the center subspace, $f$, is generally not orthogonal to the basis of the expanding subspace. That linear system is derived in Section \ref{sec:s3-derivation}. Another consequence of the three-term splitting is the emergence of the neutral contribution of linear response. Fortunately, as shown in Eq. \ref{eqn:s3-center-product}-\ref{eqn:s3-center-final}, this part of the algorithm re-uses some ingredients of the stable contribution and only requires computing $K\in\mathbb{Z}^+$ $k$-time correlations through ergodic-averaging. Finally, the evaluation of the unstable part also requires some adjustments. Eq. \ref{eqn:s3-unstable} indicates that we need $c^i$, $i=0,1,...,m$, their unstable derivatives $b$, and derivatives of the SRB measure represented by $g$. We acknowledge that the computation of the SRB measure gradient is agnostic to the presence of the center manifold. Using the measure preservation property and chain rule on smooth manifolds, one can derive exponentially converging recursive formulas for $g$. The reader is referred to the authors' previous work published in \cite{sliwiak-srb} for a detailed derivation and analysis of a trajectory-driven algorithm for $g$. Therefore, we only need to modify the way $b$ is computed in the presence of the neutral subspace. Once $b$ is found, the unstable part is computed similarly to its neutral counterpart, by summing up $K$ $k$-time correlations.

Note that $b^{i,j}$ is defined as the directional derivative of $c^i$ computed along the $j$-th basis vector $q^j$. While the regularized form of the unstable contribution (RHS of Eq. \ref{eqn:s3-unstable}) involves only self-derivatives of $c^i$, i.e., $b^{i,j}$ with $i=j$, we show that in order to find a trajectory-following recursion, we also need all possible cross-derivatives of $c^i$. The main tool used in the derivation of these formulas is the measure-based parameterization of local unstable manifolds with orthonormal gradients \cite{sliwiak-srb}. It means that the $m$-dimensional unstable manifold $U_k$ including $x_k$, i.e., the point of $M$ crossed by the trajectory at the $k$-th time step, is parameterized as follows: $x_k(\xi):[0,1]^m\to U_k\subset M$ such that $x_k(\xi)$ is the multivariate inverse cumulative distribution (quantile function) and $\nabla_{\xi_k}x_k = Q_k$. In this context, the marginal SRB density $\rho_k$ defined on $U_k$ can be viewed as the probability density function (PDF) of the uniform measure nonlinearly re-distributed by $x_k(\xi)$. The chart coordinates $\xi_k$ are updated step-by-step to ensure the orthogonality of the gradient $\nabla_{\xi_k}x_k = [\partial_{\xi_k^1}x_k,...,\partial_{\xi_k^m}x_k]$. A more rigorous description and analysis of this coordinate transformation can be found in \cite{sliwiak-srb}.  

To obtain $b^{i,j}$, $i = 0,1,...,m$, $j = 1,...,m$, we simply differentiate Eq. \ref{eqn:s3-inhomog}, Eq. \ref{eqn:s3-constraint2} and the constraint $v\cdot f = 0$ with respect to all components of $\xi$, apply the chain rule, and solve a linear system with $m(m+1)$ equations and the same number of unknowns. Notice that, assuming $\nabla_{\xi_k}x_k = Q_k$, the directional derivatives along $q^i$ are the same as parametric derivatives with respect to $\xi^i$.

Differentiation of Eq. \ref{eqn:s3-inhomog} with respect to $\xi_{k+1}^j$ yields
\begin{equation}
    \label{eqn:appa-inhomog-der}
    \begin{split}
    \partial_{\xi_{k+1}^j}v_{k+1}:&= w_{k+1}^j = \partial_{\xi_{k+1}^j}r_{k+1} - \sum_{l=1}^m b_{k+1}^{l,j}\,q_{k+1}^l + c_{k+1}^l\,p^{l,j} \\ & - b_{k+1}^{0,j}\,f_{k+1} - c_{k+1}^0\,Df_{k+1}\,q_{k+1}^j,
    \end{split}
\end{equation}
where $p^{i,j}:=\partial_{\xi^j} q^i$. In the above equation, we used the following identity, 
$$\partial_{\xi^j_{k+1}} f_{k+1} = Df_{k+1}\,\partial_{\xi_{k+1}^j} x_{k+1} = Df_{k+1}\,q_{k+1}^j.$$ 
Consequently, differentiating Eq. \ref{eqn:s3-constraint2}, i.e., constraint enforcing $v\cdot q = 0$, with respect to $\xi_{k+1}^j$ gives
\begin{equation}
    \label{eqn:appa-constraint2-der}
    \begin{split}
    b_{k+1}^{i,j} & = p_{k+1}^{i,j}\cdot\left(r_{k+1}-c_{k+1}^0\,f_{k+1}\right) + q_{k+1}^i\cdot\partial_{\xi_{k+1}^j}r_{k+1} \\ &
    - b_{k+1}^{0,j}q_{k+1}^i\cdot f_{k+1} - c_{k+1}^0\,q_{k+1}^i\cdot Df_{k+1}\,q_{k+1}^j.
    \end{split}
\end{equation}
To eliminate $w$ from the linear system, we differentiate the constraint $v\cdot f=0$ with respect to $\xi_{k+1}^j$ and plug Eq. \ref{eqn:appa-inhomog-der} to obtain
\begin{equation}
    \label{eqn:appa-inhomog-der-projection}
    \begin{split}
      -v_{k+1}\cdot Df_{k+1}\,q_{k+1}^j & = w_{k+1}^j\cdot f_{k+1} =  \partial_{\xi_{k+1}^j}r_{k+1}\cdot f_{k+1} \\ & - \sum_{l=1}^m b_{k+1}^{l,j}\,q_{k+1}^l\cdot f_{k+1} + c_{k+1}^l\,p^{l,j}\cdot f_{k+1} \\ & - b_{k+1}^{0,j}\,f_{k+1}\cdot f_{k+1} - c_{k+1}^0\,Df_{k+1}\,q_{k+1}^j\cdot f_{k+1}.
    \end{split}
\end{equation}
Finally, by combining Eq. \ref{eqn:appa-constraint2-der}--\ref{eqn:appa-inhomog-der-projection}, we derive the following linear system for $b^{i,j}$, $i = 0,1,...,m$, $j = 1,...,m$,
\begin{equation}
    \label{eqn:appa-system}
    \begin{split}
    &(f_{k+1}\cdot f_{k+1}) b_{k+1}^{0,j} + \sum_{l=1}^m (q_{k+1}^l\cdot f_{k+1})\,b^{l,j}_{k+1} = d_{k+1}^{0,j},\,\,\,j=1,...,m,\\
    &
    (q_{k+1}^i\cdot f_{k+1})b_{k+1}^{0,j} + b^{i,j}_{k+1} = d_{k+1}^{i,j},\,\,\,i,j=1,...,m,
    \end{split}
\end{equation}
where
\begin{equation}
    \label{eqn:appa-system-rhs}
    \begin{split}
    & d_{k+1}^{0,j} := v_{k+1}\cdot Df_{k+1}\,q_{k+1}^j + \partial_{\xi_{k+1}^j}r_{k+1}\cdot f_{k+1} \\ & - \sum_{l=1}^m c_{k+1}^l\,p_{k+1}^{l,j}\cdot f_{k+1} - c_{k+1}^0\,Df_{k+1}\,q_{k+1}^j\cdot f_{k+1},\,\,\,j=1,...,m, \\
    &
    d_{k+1}^{i,j} := p_{k+1}^{i,j}\cdot\left(r_{k+1}-c_{k+1}^0\,f_{k+1}\right) + q_{k+1}^i\cdot\partial_{\xi_{k+1}^j}r_{k+1} \\ 
    &
    - c_{k+1}^0\,q_{k+1}^i\cdot Df_{k+1}\,q_{k+1}^j,\,\,\,i,j=1,...,m.
    \end{split}
\end{equation}
The Schur complement of System \ref{eqn:appa-system}--\ref{eqn:appa-system-rhs} consists of $m^2$ constant-diagonal blocks. Their values are exactly the same as the corresponding entries of $S$. Therefore, if the inverse $S^{-1}$ is available, we can directly compute the sought-after quantities,
\begin{equation}
    \label{eqn:appa-system-b}
    \begin{split}
    b_{k+1}^{i,j} &= (S_{k+1}^{-1})^{i:}\cdot d_{k+1}^{1:m,j} - \sum_{l=1}^m (S_{k+1}^{-1})^{il}\frac{q_{k+1}^l\cdot f_{k+1}}{f_{k+1}\cdot f_{k+1}}d_{k+1}^{0,j}
    \\ &
    =(S_{k+1}^{-1})^{i:}\cdot\left(d_{k+1}^{1:m,j}-\frac{d_{k+1}^{0,j}}{f_{k+1}\cdot f_{k+1}}Q_{k+1}^T\,f_{k+1}\right),\,\,\,i,j=1,...,m,
    \end{split}
\end{equation}
where $(S^{-1})^{ij}$ indicates the entry of $S^{-1}$ corresponding to its $i$-th row and $j$-th column. Analogously, $d^{1:m,j}$ denotes the $m$-dimensional array including all $d^{i,j}$ for all $i=1,...,m$ and a fixed $j$. Once $b^{i,j}$ for all $i,j=1,...,m$ is computed, $b^{0,j}$ and $w^j$, $j=1,...,m$ can be evaluated directly using Eq. \ref{eqn:appa-inhomog-der} and Eq. \ref{eqn:appa-system}.

Based on Eq. \ref{eqn:appa-inhomog-der}--\ref{eqn:appa-system-b}, we can now construct a trajectory-following iteration to compute $b$. These equations involve some ingredients previously derived for the stable and neutral parts. The new quantities are the parametric derivatives of the basis vectors $p$, i.e., derivative of Lyapunov vectors, and $\partial_{\xi_{k+1}^j} r_{k+1}$. The former are computed using the procedure for $g$ extended by an extra low-cost projection \cite{sliwiak-s3}. Using the definition of $r_{k+1}$ and all underlying quantities, we apply the chain rule to expand $\partial_{\xi_{k+1}^j} r_{k+1}$,
\begin{equation}
    \label{eqn:appa-chain-rule}
    \partial_{\xi_{k}^j} r_{k+1} = D^2\varphi(v_k,q_k^j) + D\varphi_k\,w_k^j + D\partial_s\varphi_k\,q_k^j,
\end{equation}
where $D^2\varphi(a,b)$ denotes the second-order bilinear form whose $i$-th component equals $(D^2\varphi(a,b))^i = \partial_{x^k}\partial_{x^l} \varphi^i\,a^k\,b^l$ (per Einstein's summation convention), while $D\partial_s \varphi$ denotes the phase-space Jacobian of parametric derivative of $\varphi$. Note also that Eq. \ref{eqn:appa-chain-rule} needs to be further re-scaled by the Jacobian of the coordinate transformation from $\xi_{k}$ to $\xi_{k+1}$. Without loss of generality, we can choose $\xi = 0$ and show that the Jacobian of coordinate transformation is a by-product of the iterative algorithm for the basis vectors $q$ \cite{sliwiak-srb}. Based on the above derivations, Section \ref{sec:s3-derivation} and \cite{sliwiak-s3}, Algorithm \ref{alg:alg1} summarizes all the steps required to approximate the full linear response of a hyperbolic flow. While the most important aspects are covered in this work, the reader is referred to these two external references for a rigorous justification of all other parts. 

\begin{algorithm}\label{alg:alg1}
\SetAlgoLined
\SetKwInOut{Input}{Input}
\SetKwInOut{Output}{Output}
\Input{$N$, $K$, $T$, $n$, $m$}
\Output{$d\langle J\rangle/ds \approx (s + c + u)/N$}

Randomly generate: $x_0$, $v_0$, $Q_0$, $a^{i,j}_0$, $w^i_0$ for all $i,j=1,...,m$\;
Set $s=c=u=0$\;

\For(\tcp*[h]{main time loop}){$k = 0,...,N-1$}{

\If{$k\geq T$}
{$s := s + DJ_k\cdot v_k$\;
$u := u - J_k\,(u_k + u_{k-1} + ... + u_{k-K+1});$
        \hspace*{1.2em}%
        \rlap{\smash{$\left.\begin{array}{@{}c@{}}\\{}\\{}\\{}\\{}\end{array}\color{red}\right\}%
          \color{red}\begin{tabular}{l}Update stable ($s$), neutral ($c$)\\and unstable ($u$) contributions \cite{chandramoorthy-phdthesis,sliwiak-s3}.\end{tabular}$}}\\
$c := c + DJ_k\cdot f_{k}\,(c_k^0 + c_{k-1}^0 + ... + c_{k-K+1}^0)$\;}
$P_{k+1} = D\varphi_k\,Q_k$\;
QR-factorize $P_{k+1}$: $Q_{k+1}\,R_{k+1} = P_{k+1};$ 
\hspace*{-.7em}%
        \rlap{\smash{$\left.\begin{array}{@{}c@{}}\\{}\\{}\end{array}\color{red}\right\}%
          \color{red}\begin{tabular}{l}Update basis vectors $Q$ and transformation\\Jacobian $R$. See \cite{ershov-lyapunov} for derivation \\ and convergence analysis of Lyapunov basis.\end{tabular}$}}\\
Find the inverse of $R_{k+1}$\;

\lFor{$i = 1,...,m$, $j= 1,...,i$}{
$\tilde{a}_{k+1}^{i,j} = D^2\varphi_k(q_k^{i}, q_k^{j}) + D\varphi_k\,a_{k}^{i,j}$
}

\lFor{$i = 1,...,m$, $j = 1,...,i$}{
$a_{k+1}^{i,j} = \tilde{a}_{k+1}^{p,q}\,(R^{-1})_{k+1}^{pi}\,(R^{-1})_{k+1}^{qj}$
}

\For{$i = 1,...,m$}{
\For{$p, q = 1,...,m$}{
$(\partial_{\xi_{k+1}^i}R_{k+1})^{pq} = \begin{cases} q_{k+1}^p\cdot a_{k+1}^{p,i}, & \text{if } p = q 
\hspace*{1em}%
        \rlap{\smash{$\left.\begin{array}{@{}c@{}}\\{}\\{}\\{}\\{}\\{}\\{}\\{}\\{}\end{array}\color{red}\right\}%
          \color{red}\begin{tabular}{l}Compute all $m$ \\ components of $g$.\\ This part requires \\ solving $m^2$ $2^{nd}$-order\\ tangent equations.\\ Run once per time step, \\regardless of $\mathrm{dim}(J)$ \\ and $\mathrm{dim}(s)$. See \cite{sliwiak-srb} \\ for original derivation.\end{tabular}$}}
\\ q_{k+1}^p\cdot a_{k+1}^{q,i} + q_{k+1}^{q}\cdot a_{k+1}^{p,i}, & \text{if } p < q \\
    0, & \text{otherwise}\end{cases}
$\;
}
$g_{k+1}^i = -\text{tr}(\partial_{\xi_{k+1}^i}R_{k+1})$\;
}

\lFor{$i,j = 1,...,m$}{
$p_{k+1}^{i,j} = a_{k+1}^{i,j} - q_{k+1}^{l}(\partial_{\xi_{k+1}^j}R_{k+1})^{li}$ \hspace*{-0.2em}%
        \rlap{\smash{$\left.\begin{array}{@{}c@{}}\\{}\end{array}\color{red}\right\}%
          \color{red}\begin{tabular}{l}Compute derivatives of \\ Lyapunov vectors \cite{sliwiak-s3}.\end{tabular}$}}}

$S_{k+1} = I - Q_{k+1}^T f_{k+1}(Q_{k+1}^T f_{k+1})^T / f_{k+1}\cdot f_{k+1}$\;

Find the inverse of $S_{k+1}$\;

$r_{k+1} = D\varphi_k\,v_{k} + \chi_{k+1}$\;

$z_{k+1} = Q_{k+1}^T\left(r_{k+1} - (f_{k+1}\cdot r_{k+1})/(f_{k+1}\cdot f_{k+1})f_{k+1}\right)$\;

\lFor{$i = 1,...,m$}{
$c_{k+1}^i = (S_{k+1}^{-1})^{ij}\,z_{k+1}^j$}

$c_{k+1}^0 = f_{k+1}\cdot(r_{k+1}-c_{k+1}^i\,q_{k+1}^i)/f_{k+1}\cdot f_{k+1}$\;

$v_{k+1} = r_{k+1} - c_{k+1}^i\,q_{k+1}^i - c_{k+1}^0\,f_{k+1}$\;

\lFor{$i = 1,...,m$}{
$\partial_{\xi_k^i}\,r_{k+1} = D^2\varphi_k(v_k,q_k^i) + D\varphi_k\,w_k^i + D\partial_s\varphi_k\,q_k^{i}$}

$\nabla_{\xi_{k+1}}r_{k+1} = \nabla_{\xi_{k}}r_{k+1}\,R_{k+1}^{-1};$ \hspace*{34.2em}%
        \rlap{\smash{$\left.\begin{array}{@{}c@{}}\\{}\\{}\\{}\\{}\\{}\\{}\\{}\\{}\\{}\\{}\\{}\\{}\\{}\\{}\\{}\\{}\\{}\\{}\\{}\\{}\\{}\end{array}\color{red}\right\}%
          \color{red}\begin{tabular}{l}Modification \\ of the general \\ discrete S3 \\from \cite{sliwiak-s3}.\\ This part \\ computes \\ $c$, $b$, $w$,\\according to \\ the derivation\\presented in\\ this section \\and Section 2.1.\end{tabular}$}}\\

\lFor{$i = 1,...,m$}{
$d_{k+1}^{0,i} = v_{k+1}\cdot Df_{k+1}\,q_{k+1}^i + \partial_{\xi_{k+1}^i}r_{k+1}\cdot f_{k+1} - c_{k+1}^l\,p_{k+1}^{l,i}\cdot f_{k+1} - c_{k+1}^0\,Df_{k+1}\,q_{k+1}^i\cdot f_{k+1}$}

\lFor{$i,j = 1,...,m$}{
$d_{k+1}^{i,j} = p_{k+1}^{i,j}\cdot (r_{k+1}-c_{k+1}^0\,f_{k+1}) + q_{k+1}^i\cdot\partial_{\xi_{k+1}^j}r_{k+1}
    - c_{k+1}^0\,q_{k+1}^i\cdot Df_{k+1}\,q_{k+1}^j$}

\lFor{$i,j = 1,...,m$}{
$b_{k+1}^{i,j} = (S_{k+1}^{-1})^{i:}\cdot (d_{k+1}^{1:m,j}- d_{k+1}^{0,j}/(f_{k+1}\cdot f_{k+1})Q_{k+1}^T\,f_{k+1})$}

\lFor{$i = 1,...,m$}{
$b_{k+1}^{0,i} = 1/(f_{k+1}\cdot f_{k+1})\,(d_{k+1}^{0,i} - (q_{k+1}^l\cdot f_{k+1})b_{k+1}^{l,i})$}

\lFor{$i = 1,...,m$}{
$w_{k+1}^{i} = \partial_{\xi_{k+1}^i} r_k - b_{k+1}^{l,i}\,q_{k+1}^l - c_{k+1}^l\,p_{k+1}^{l,i} - b_{k+1}^{0,i}\,f_{k+1} - c_{k+1}^0\,Df_{k+1}q_{k+1}^i$
}

Save the two scalars: $u_{k+1} = b^{i,i}_{k+1} + c_{k+1}^i\,g_{k+1}^i$ and $c_{k+1}^0$\;
Advance the iteration: $x_{k+1} = \varphi(x_k)$\;
}

\caption{Space-split sensitivity algorithm for hyperbolic flows}
\end{algorithm}

The input parameter $T$ is to allow all the recursions to converge before the linear response contributions are collected. Note that Algorithm \ref{alg:alg1} is agnostic to the time integration method, which directly affects $\varphi$ and hence the cost of computing its derivatives. In \ref{sec:derivatives}, we derive relevant differentiation operators for the midpoint scheme.

Assuming both the objective function $J$ and parameter $s$ are scalars, the computational cost of Algorithm \ref{alg:alg1} depends on three parameters: the trajectory length $N$, dimension of both the system $n$ and unstable subspace $m$. In this case, the most expensive part is the computation of the SRB density gradient (Lines 12-18). This chunk of the algorithm solves $m^2$ second-order tangent equations (Line 12) and performs double contraction against the transformation Jacobian (Line 13) to stabilize the iteration, which costs $\mathcal{O}(n^3\,m^2 + n\,m^3)$ floating point operations (flops) per time step. If $s$ is an $n_s$-dimensional vector, then the majority of the modified part of Algorithm \ref{alg:alg1} (Lines 22-35) will need to be repeated $n_s$ times, which costs $\mathcal{O}(n_s\,(n^3\,m + m^2\,n))$ flops per time step. Finally, Lines 4-8 would need to be repeated $n_J$ times if $J$ was an $n_J$-dimensional vector. This would incur an extra cost proportional to $\mathcal{O}(n_J\,n_s\,n)$ flops. Therefore, assuming $\max(m,n_s,n_J)\ll n$, the leading flop count term of the total cost of Algorithm \ref{alg:alg1} is
%\begin{equation}
%    \label{eqn:alg-cost}
%    \mathcal{O}(n^3\,m^2 + n\,m^3) + s\,\mathcal{O}(n^3\,m + m^2\,n) + t\,\mathcal{O}(n\,s).
%\end{equation}
\begin{equation}
    \label{eqn:appa-cost}
    \mathcal{O}\left(N\,n^3\,(m^2 + n_s\,m)\right). 
\end{equation}
Note that the most important factor in determining the total cost is the system's dimension $n$. This number is cubed because of the contraction of the second-order operator with two different vectors (Line 12). In practice, however, the linear differentiation operators (Jacobians, Hessians) have sparse/banded structure. This usually happens in case of PDE-related dynamical systems that have been derived using standard discretization methods such as the finite element method. The major consequence of the local structure is that the cost of evaluating first- and second-order operator-vector contractions is in fact linear to the dimension of the system. Therefore, the leading term of the flop count dramatically decreases to 
\begin{equation}
    \label{eqn:appa-cost2}
    \mathcal{O}\left(N\,n\,(m^3 + n_s\,m^2 + n_s\,n_J)\right).
\end{equation}

\section{Computing derivatives of $\varphi$ and implicit schemes}\label{sec:derivatives}

Both Algorithm \ref{alg:alg_reduced} and Algorithm \ref{alg:alg1} require computing first-order derivatives in phase space as well as parametric derivatives of $\varphi$. The latter also requires second-order derivatives to compute $g$ and $w$. They are products of the chain rule applied to the discrete version of the time-continuous system. The computational cost of evaluating these quantities heavily depends on the time integrator. For the Euler method, for example, differentiation of $\varphi$ is equally expansive as differentiation of $f$. In this paper, we use second- and fourth-order fully-explicit Runge-Kutta schemes, which involve nested functions. If the system is sparse and its dimension $n$ is large, it is efficient to compute all the tensor-vector contractions as we go rather than evaluating and storing large Jacobians and Hessians. Therefore, our aim is to use the chain rule to express all contraction types appearing in both algorithms such as $D\varphi\,v$ in terms of similar tensor-vector products involving derivatives of $f$ only. In this section, we present derivations for the second-order Runge-Kutta map defined by Eq. \ref{eqn:s3-time-int}. Analogous expressions for the fourth-order scheme can be found in the attached Python code.

For the midpoint method, $\varphi(x_k)$ is defined as
\begin{equation}
    \label{eqn:appb-map}
    \varphi(x_k) = x_k + \Delta t\,f\left(x_k + \frac{\Delta t}{2}\,f(x_k)\right):=x_k + \Delta t\,f(x_p) = x_{k} + \Delta t\,f_p,
\end{equation}
where $x_p:=x_k + \Delta t/2\,f(x_k)$. Therefore, for any vector $v\in \mathbb{R}^n$, 
\begin{equation}
    \label{eqn:appb-dphi}
    D\varphi_k\,v = v + \Delta t\,Df_p\,v + \frac{\Delta t^2}{2}\,Df_p\,Df_k\,v,
\end{equation}
with $Df_k = Df(x_k)$ and $Df_p = Df(x_p)$, in consistency with our notation convention. Differentiating Eq. \ref{eqn:appb-dphi} once more and contracting it against yet another vector $a\in \mathbb{R}^n$, we obtain the following relation, 
\begin{equation}
    \label{eqn:appb-d2phi}
    \begin{split}
    D^2\varphi_k(v,a) =& \Delta t\,D^2 f_p\left(v + \frac{\Delta t}{2}\,Df_k\,v,a\right) + \\& \frac{\Delta t^2}{2}\,D^2 f_p\left(v + \frac{\Delta t}{2}\,Df_k\,v,Df_k\,a\right) + \frac{\Delta t^2}{2}\,Df_p\,D^2f_k(v,a).
    \end{split}
\end{equation}
Recall that $D^2\varphi_k(v,a) \in \mathbb{R}^n$. Assuming $f$ also depends on a scalar parameter $s$, the parametric derivative of Eq. \ref{eqn:appb-map} expands as follows,
\begin{equation}
    \label{eqn:appb-dsphi}
    \partial_s\varphi_k = \chi_{k+1} = \Delta t\,\partial_s f_p + \frac{\Delta t^2}{2}Df_p\,\partial_s f_k,
\end{equation}
where $\partial_s f_k = \partial f/\partial s\,(x_k)$. The final relevant contraction, $D\partial_s\varphi_k\,v$, involves mixed parametric and phase-space derivatives and is obtained by differentiating Eq. \ref{eqn:appb-dsphi}, 
\begin{equation}
    \label{eqn:appb-ddsphi}
    \begin{split}
    D\partial_s\varphi_k\,v =& \Delta t\,D\partial_s f_p \left(v + \frac{\Delta t}{2}\,Df_k\,v\right) + \\&\frac{\Delta t^2}{2}D^2f_p\left(v + \frac{\Delta t}{2}\,Df_k\,v ,\partial_s f_k\right) + \frac{\Delta t^2}{2}\,Df_p\,D\partial_s f\,v.
    \end{split}
\end{equation}

%Summarize (max. two sentences) and comment on implicit schemes
We highlight the fact that, for the midpoint method, each tensor-vector product involving $\varphi$ requires the evaluation of $\mathcal{O}(1)$ similar products containing $f$. The fourth-order Runge-Kutta scheme is in fact a four-level nested map from $x_k$ to $x_{k+1}$. In this case, the Hessian-vectors contraction requires about 20 such evaluations. For sparse systems, however, the cost of a single evaluation of $Df\,v$, $D^2f (a,v)$, $D\partial_s f\,v$ is linear in $n$.  

An implicit scheme is a common choice for stiff systems. That choice does not affect our linear response algorithms. The only part that needs to be modified is the way the products appearing in Eq. \ref{eqn:appb-map}--\ref{eqn:appb-ddsphi} are computed. Let us consider a generic implicit scheme,
\begin{equation}
    \label{appb-implicit}
    h(x_k,x_{k+1}) = 0,
\end{equation}
where $x_{k+1} = \varphi(x_k)$. Assuming $x_k$ is known, the $n$-dimensional nonlinear system defined by Eq. \ref{appb-implicit} is typically solved for $x_{k+1}$ using a standard solver such as the Newton-Raphson method. Differentiating Eq. \ref{appb-implicit} with respect to $x_{k}$ and multiplying both sides by a vector $v$, we obtain the following system,
\begin{equation}
    \label{appb-system}
    \frac{\partial h}{\partial x_{k+1}}\,D\varphi_k\,v = -\frac{\partial h}{\partial x_{k}}\,v,
\end{equation}
where $\partial h /\partial x_{k}$ and $\partial h /\partial x_{k+1}$ respectively the $n\times n$ Jacobian matrices of $h$ with respect to $x_{k}$ and $x_{k+1}$, respectively, both evaluated at $(x_k,x_{k+1})$. If both $x_{k}$ and $x_{k+1}$ are known, the linear system defined by Eq. \ref{appb-system} can be solved for $D\varphi_k\,v$, which is a necessary ingredient of our linear response algorithms. To compute other tensor-vector products, we further differentiate Eq. \ref{appb-system}, apply the chain rule as presented above, and formulate analogous linear systems.   

%-------------------------End of appendices----------------------

\end{document}